\documentclass[10pt,reqno]{amsart}
\pdfoutput=1 
\usepackage[a4paper,margin=3cm]{geometry}
\usepackage[T1]{fontenc}
\usepackage[utf8]{inputenc}
\usepackage[pagebackref]{hyperref}\hypersetup{pdfborder=0 0 0}
\usepackage{latexsym,amssymb,lmodern,graphicx,caption,listings}
\usepackage[francais,english]{babel}

\title{From Boltzmann to random matrices and beyond}

\author{Djalil~Chafaï} %
\address{CEREMADE, UMR 7534, Université Paris Dauphine, PSL,
  IUF, France.} %
\email{\url{mailto:djalil(at)chafai.net}} %
\urladdr{\url{http://djalil.chafai.net/}}

\date{Spring 2014. Revised July 2, October 31, for Ann. Fac. Sci. Toulouse. Compiled
  \today.}

\newtheorem{theorem}{Theorem}[section]%
\newtheorem{remark}[theorem]{Remark}%
\newcommand{\dC}{\mathbb{C}}
\newcommand{\dE}{\mathbb{E}}

\newcommand{\dN}{\mathbb{N}}
\newcommand{\dP}{\mathbb{P}}
\newcommand{\dR}{\mathbb{R}}

\newcommand{\dZ}{\mathbb{Z}}

\newcommand{\cA}{\mathcal{A}}\newcommand{\cB}{\mathcal{B}}
\newcommand{\cC}{\mathcal{C}}
\newcommand{\cF}{\mathcal{F}}
\newcommand{\cG}{\mathcal{G}}\newcommand{\cH}{\mathcal{H}}
\newcommand{\cI}{\mathcal{I}}
\newcommand{\cK}{\mathcal{K}}
\newcommand{\cM}{\mathcal{M}}\newcommand{\cN}{\mathcal{N}}

\newcommand{\cR}{\mathcal{R}}
\newcommand{\cS}{\mathcal{S}}

\newcommand{\al}{\alpha}
\newcommand{\be}{\beta}
\newcommand{\De}{\Delta}
\newcommand{\de}{\delta}

\newcommand{\Ga}{\Gamma}
\newcommand{\la}{\lambda}

\newcommand{\si}{\sigma}

\newcommand{\beq}{\begin{equation}}
\newcommand{\eeq}{\end{equation}}
\newcommand{\brem}{\begin{remark}}
\newcommand{\erem}{\end{remark}}

\newcommand{\veps}{\varepsilon}

\newcommand{\varn}{\varnothing}
\newcommand{\ABS}[1]{{{\left| #1 \right|}}} 
\newcommand{\BRA}[1]{{{\left\{#1\right\}}}} 
\newcommand{\DOT}[1]{{{\left<#1\right>}}} 
\newcommand{\NRM}[1]{{{\left\| #1\right\|}}} 
\newcommand{\PAR}[1]{{{\left(#1\right)}}} 
\newcommand{\pd}{{\partial}} 
\newcommand{\SBRA}[1]{{{\left[#1\right]}}} 
\newcommand{\LIP}[1]{\|#1\|_{\mathrm{Lip}}} 

\newcommand{\SUPP}{\mathrm{supp}}

\newcommand{\TR}{\mathrm{Tr}}

\newcommand{\IND}{\mathbf{1}}

\newcommand{\OL}[1]{\overline{#1}}

\renewcommand{\Im}{\mathfrak{Im}}
\renewcommand{\Re}{\mathfrak{Re}}

\numberwithin{equation}{section}



\begin{document}

\begin{abstract}
  These expository notes propose to follow, across fields, some aspects of the
  concept of entropy. Starting from the work of Boltzmann in the kinetic
  theory of gases, various universes are visited, including Markov processes
  and their Helmholtz free energy, the Shannon monotonicity problem in the
  central limit theorem, the Voiculescu free probability theory and the free
  central limit theorem, random walks on regular trees, the circular law for
  the complex Ginibre ensemble of random matrices, and finally the asymptotic
  analysis of mean-field particle systems in arbitrary dimension, confined by
  an external field and experiencing singular pair repulsion. The text is
  written in an informal style driven by energy and entropy. It aims to be
  recreative and to provide to the curious readers entry points in the
  literature, and connections across boundaries.
  \medskip

  \begingroup\selectlanguage{francais}
  \textsc{Résumé.} Ces notes d'exposition proposent de suivre, à travers
  différents domaines, quelques aspects du concept d'entropie. À partir du
  travail de Boltzmann en théorie cinétique des gas, plusieurs univers sont
  visités, incluant les processus de Markov et leur énergie libre de
  Helmholtz, le problème de Shannon de monotonie de l'entropie dans le
  théorème central limite, la théorie des probabilités libres de Voiculescu et
  le théorème central limite libre, les marches aléatoires sur les arbres
  réguliers, la loi du cercle pour l'ensemble de Ginibre complexe de matrices
  aléatoires, et enfin l'analyse asymptotique de systèmes de particules champ
  moyen en dimension arbitraire, confinées par un champ extérieur et subissant
  une répulsion singulière à deux corps. Le texte est écrit dans un style
  informel piloté par l'énergie et l'entropie. Il vise a être récréatif, à
  fournir aux lecteurs curieux des points d'entrée dans la littérature, et des
  connexions au delà des frontières.
  \endgroup\selectlanguage{english}
\end{abstract}


\maketitle

{\small\tableofcontents}

This text forms the written notes of a talk entitled ``About confined
particles with singular pair repulsion'', given at the occasion of the
workshop ``Talking Across Fields'' on convergence to the equilibrium of Markov
chains. This workshop was organized in Toulouse from 24 to 28 March 2014 by
Laurent Miclo, at the occasion of the CIMI Excellence research chair for Persi
Diaconis.

Almost ten years ago, we wanted to understand by curiosity the typical global
shape of the spectrum of Markov transition matrices chosen at random in the
polytope of such matrices. It took us several years to make some progress
\cite{MR2549497,MR2892961,MR3168123}, in connection with the circular law
phenomenon of Girko. The circular law states that the empirical measure of the
eigenvalues of a random $n\times n$ matrix, with i.i.d.\ entries\footnote{Real
  or complex, it does not matter, with an arbitrary mean, as soon as we look
  at narrow convergence.} of variance $1/n$, tends to the uniform law on the
unit disc of the complex plane, as the dimension $n$ tends to infinity. This
universal result was proved rigorously by Tao and Vu \cite{MR2722794}, after
fifty years of contributions. The proof of this high dimensional phenomenon
involves tools from potential theory, from additive combinatorics, and from
asymptotic geometric analysis. The circular law phenomenon can be checked in
the Gaussian case using the fact that the model is then exactly solvable.
Actually, Ginibre has shown in the 1960's that if the entries are i.i.d.\
centered complex Gaussians then the eigenvalues of such matrices form a
Coulomb gas at temperature $1/n$ in dimension $2$. This in turn suggests to
explore the analogue of the circular law phenomenon in dimension $\geq3$,
beyond random matrices. This leads us to introduce in
\cite{2013arXiv1304.7569C} stochastic interacting particle systems in which
each particle is confined by an external field, and each pair of particles is
subject to a singular repulsion. Under general assumptions and suitable
scaling, the empirical measure of the particles converges, as the number of
particles tends to infinity, to a probability measure that minimizes a natural
energy-entropy functional. In the case of quadratic confinement and Coulomb
repulsion, the limiting law is uniform on a ball.

This expository text is divided into five sections, written in an informal
style. The first section introduces the Boltzmann entropy and H-Theorem, and
the analogous concept of Helmholtz free energy for Markov processes. The
second section discusses the Shannon monotonicity problem of the Boltzmann
entropy along the central limit theorem. The third section introduces
Voiculescu free probability theory, the free entropy, and the free central
limit theorem. The fourth section discusses the circular law for random
matrices drawn from the complex Ginibre ensemble, using a large deviations
rate function inspired from the Voiculescu entropy. The fifth and last section
goes beyond random matrices and studies mean-field particle systems with
singular interaction, for which the large deviations rate function is
connected to Boltzmann and Voiculescu entropies.

Talking across fields is rarely an easy task. You may know for instance that
Andreï Andreïevitch Markov (1856--1922) published his seminal article on what
we call Markov chains in 1906\footnote{At the age of fifty, the year of the
  death of Ludwig Boltzmann.}, and that approximately at the same time, the
theory of non-negative matrices was developed by Oskar Perron (1880--1975) and
Ferdinand Georg Frobenius (1849--1917). It took several years to talk across
fields, and according to Eugene Seneta \cite{seneta-markov-and-the-birth}, the
link was made by von Mises (1883--1953). Various point of views are available
on Markov chains. Beyond the concreteness of conditional construction and
stochastic simulation, a Markov model can always be seen as a deterministic
evolution, along the time, of a probability distribution. This mechanical view
of (random) nature can be traced back to Charles Darwin (1809--1882) with his
mutation-selection evolution theory, and to Ludwig Boltzmann (1844--1906) with
his H-Theorem in atomistic kinetic theory of gases. These two great figures of
the nineteenth century can be viewed as proto-Markovian.

\medskip

\begin{quote}
  ``\emph{If you ask me about my innermost conviction whether our century will
    be called the century of iron or the century of steam or electricity, I
    answer without hesitation: It will be called the century of the mechanical
    view of Nature, the century of Darwin.}''
  \begin{flushright}
    Ludwig Boltzmann, Vienna, 1886\\
    Der zweite Hauptsatz der mechanischen Wärmetheorie\\
    Lecture at the Festive Session, Imperial Academy of Sciences
  \end{flushright}
\end{quote}

\medskip

\begin{center}
  \begin{figure}[htbp]
    \includegraphics[scale=0.5]{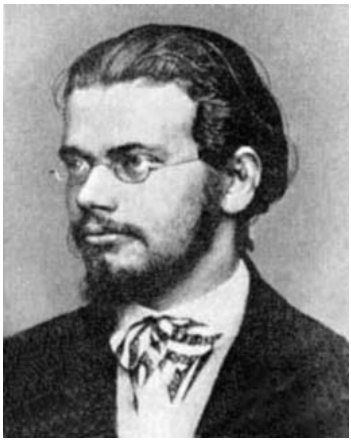}
    \caption*{Ludwig Boltzmann (1844--1906)}
  \end{figure}
\end{center}

\begin{quote}
  ``\emph{Boltzmann summarized most (but not all) of his work in a two volume
    treatise Vorlesungen über Gastheorie. This is one of the greatest books in
    the history of exact sciences and the reader is strongly advised to
    consult it. It is tough going but the rewards are
    great\footnote{Personally, when I was a first year graduate student at
      Université Paul Sabatier in Toulouse, I borrowed this book from the
      university library by curiosity. I can confirm the ``tough going''
      aspect mentioned by Mark Kac!}.}'' %
  \begin{flushright}
    Mark Kac, Ithaca, 1959, excerpt from \cite{MR0102849}
  \end{flushright}
\end{quote}

\section{Ludwig Boltzmann and his H-Theorem}

\smallskip


\subsection{Entropy}

A simple way to introduce the Boltzmann entropy is to use the concept of
combinatorial disorder. More precisely, let us consider a system of $n$
distinguishable particles, each of them being in one of the $r$ possible
states (typically energy levels). We have $n=n_1+\cdots+n_r$ where $n_i$ is
the number of particles in state $i$. The vector $(n_1,\ldots,n_r)$ encodes
the macroscopic state of the system, while the microscopic state of the system
is encoded by the vector $(b_1,\ldots,b_n)\in\{1,\ldots,r\}^n$ where $b_i$ is
the state of the $i$-th particle. The number of microscopic states compatible
with a fixed macroscopic state $(n_1,\ldots,n_r)$ is given by the multinomial
coefficient\footnote{Encoding the occurrence of each face of an $r$-faces dice
  thrown $n$ times.} $n!/(n_1!\cdots n_r!)$. This integer measures the
microscopic degree of freedom given the macroscopic state. As a consequence,
the additive degree of freedom per particle is then naturally given by
$(1/n)\log(n!/(n_1!\cdots n_r!))$. But how does this behave when $n$ is large?
Let us suppose simply that $n$ tends to $\infty$ while $n_i/n\to p_i$ for
every $1\leq i\leq r$. Then, thanks to the Stirling formula, we get, denoting
$p:=(p_1,\ldots,p_r)$,
\[
\cS(p)
:=\lim_{n\to\infty}\frac{1}{n}\log\PAR{\frac{n!}{n_1!\cdots n_r!}}
=-\sum_{i=1}^rp_i\log(p_i).
\]
The quantity $\cS(p)$ is the Boltzmann entropy of the discrete probability
distribution $p$. It appears here as an asymptotic additive degree of freedom
per particle in a system with an infinite number of particles, each of them
being in one of the $r$ possible states, with population frequencies
$p_1,\ldots,p_r$. This is nothing else but the first order asymptotic analysis
of the multinomial combinatorics:
\[
\frac{n!}{n_1!\cdots n_r!}\approx e^{n\cS(n_1/n,\ldots,n_r/n)}.
\]
When the disorder of the system is better described by a probability density
function $f:\dR^d\to\dR_+$ instead of a discrete probability measure, we may
introduce by analogy, or passage to the limit, the continuous Boltzmann
entropy of $f$, denoted $\cS(f)$, or $-H(f)$ in the terminology of Boltzmann,
\[
\cS(f):=-\int_{\dR^d}\!f(x)\log(f(x))\,dx.
\]
When $X$ is a random variable, we denote by $\cS(X)$ the entropy of its law.
Here we use the notation $\cS$, which is initially the one used by Clausius for
the concept of entropy in thermodynamics.

\subsection{Maximum entropy under constraints}

The Boltzmann entropy $\cS$ measures an average disorder. One can seek for a
probability density $f_*$ that maximizes the linear functional $f\mapsto \cS(f)$
over a convex class $\cC$ formed by a set of constraints on $f$:
\[
\cS(f_*)=\max\{\cS(f):f\in\cC\}.
\] 
The class $\cC$ is typically defined by linear (in $f$) statistics on the form
$\displaystyle\int\!g(x)f(x)\,dx=c_g$ for $g\in\cG$.

Following Boltzmann, suppose that the internal state of an isolated system is
described by a parameter $x\in\dR^d$ which is statistically distributed
according to a probability density $f$, and suppose furthermore that the
energy of state $x$ is $V(x)$. Then the average energy of the system is
\[
a=\int\!V(x)\,f(x)\,dx.
\]
Let $\cC$ be the class of probability densities $f$ which satisfy this
constraint, and let us seek for $f_*\in\cC$ that maximizes the entropy $\cS$
on $\cC$, in other words such that $\cS(f_*)=\max_{f\in\cC}\cS(f)$. A Lagrange
variational analysis leads to $-\log f_*=\al+\be V$ where $\al,\be$ are
Lagrange multipliers. We select $\al,\be>0$ in such a way that $f_*\in\cC$,
which gives a unique solution
\[
f_*=\frac{1}{Z_\be}e^{-\be V}\quad\text{where}\quad Z_\be:=\int\!e^{-\be V(x)}\,dx.
\]
In Physics $\be$ is interpreted as an inverse temperature times a universal
constant called the Boltzmann constant, selected in such a way that
$f_*\in\cC$. Indeed, by using the definition of $f_*$, the fact $f,f_*\in\cC$,
and the Jensen inequality for the convex function $u\geq0\mapsto u\log(u)$, we
have
\[
\cS(f_*)-\cS(f)
=\int\!\frac{f}{f_*}\log\PAR{\frac{f}{f_*}}f_*\,dx
\geq0.
\]
The quantity in the middle is known as the Kullback-Leibler divergence or
relative entropy\footnote{This concept was actually introduced by Solomon
  Kullback (1907 -- 1994) and Richard Leibler (1914 -- 2003) in the 1950's as
  an information gain, and was inspired from the entropy in the Shannon theory
  of communication.} with respect to $f_*dx$, see \cite{MR39968,MR1461541}.
The Jensen inequality and the strict convexity of $u\geq0\mapsto u\log(u)$
give that $f_*$ is the unique density which achieves $\max_{\cC}\cS$. We write
$f_*=\arg\max_{\cC}\cS$. The lower is the energy $V(x)$ of state $x$, the
higher is the value $f_*(x)$ of the maximum entropy density $f_*$. Taking for
instance $V(x)=+\infty\IND_{K^c}(x)$ reveals that uniform laws maximize
entropy under support constraint, while taking $V(x)=\NRM{x}_2^2$ reveals that
Gaussian laws maximize entropy under second moment constraint. In particular,
on $\dR$, denoting $G$ a Gaussian random variable,
\[
\dE(X^2)=\dE(G^2)
\Rightarrow
\cS(X)\leq \cS(G)
\quad\text{and}\quad
\cS(X)=\cS(G)
\Rightarrow
X\overset{d}{=}G.
\]
It is well known in Bayesian statistics that many other classical discrete or
continuous laws are actually maximum entropy laws over classes of laws defined
by natural constraints.

\subsection{Free energy and the law of large numbers}
Still with $f_*=Z_\be^{-1}e^{-\be V}$, we have
\[
-\frac{1}{\be}\log(Z_\be)
=\cA(f_*)
\quad\text{where}\quad
\cA(f):=\int\!V(x)f(x)\,dx-\frac{1}{\be}\cS(f).
\]
The functional $\cA$ is the Helmholtz\footnote{Hermann Ludwig Ferdinand von
  Helmholtz (1821 -- 1894), inventor, among other things, of the unified
  concept of energy and its conservation in physics, in competition with
  Julius Robert von Mayer (1814 -- 1878).} free energy\footnote{Should not be
  confused with Gibbs free energy (free enthalpy) even if they are closely
  related for ideal gases.}: mean energy minus temperature times entropy. The
functional $\cA$ is essentially $-\cS$ penalized by the average energy. Also,
the functional $\cA$ admits $f_*$ as a unique minimizer over the class of
densities, without constraints. Indeed, the Helmholtz free energy is connected
to the Kullback-Leibler relative entropy: for any density $f$,
\[
\cA(f)-\cA(f_*)
=\frac{1}{\be}\int\!\frac{f}{f_*}\log\PAR{\frac{f}{f_*}}f_*\,dx\geq0
\]
with equality if and only if $f=f_*$ thanks to the strict convexity of
$u\mapsto u\log(u)$. When $f$ and $f_*$ have same average energy, then we
recover the formula for the Boltzmann entropy. As we will see later on, the
Helmholtz free energy $\cA$ plays a role for Markov processes. It emerges also
from the strong law of large numbers. More precisely, let us equip the set
$\cM_1$ of probability measures on $\dR^d$ with the narrow topology, which is
the dual topology with respect to continuous and bounded functions.
If $X_1,\ldots,X_N$ are i.i.d.\ random variables with 
law $\mu_*\in\cM_1$, then their empirical distribution
$\mu_N:=\frac{1}{N}\sum_{k=1}^n\de_{X_k}$ is a random variable on $\cM_1$, and
an asymptotic analysis due to Ivan Sanov (1919 -- 1968) in the 1950's reveals
that for every Borel set $A\subset\cM_1$, as $N\gg1$,
\[
\dP(\mu_N\in A)\approx \exp\PAR{-N\inf_{A}\cK}
\quad\text{where}\quad
\cK(\mu):=\int\!\frac{d\mu}{d\mu_*}\log\PAR{\frac{d\mu}{d\mu_*}}\,d\mu_*.
\]
The rigorous version, known as the Sanov theorem, says more precisely (see \cite{MR2571413}
for a proof) that
\[
-\inf_{\mathrm{int}(A)}\cK
\leq
\liminf_{N\to\infty}\frac{\log\dP(\mu_N\in A)}{N}
\leq
\limsup_{N\to\infty}\frac{\log\dP(\mu_N\in A)}{N}
\leq
-\inf_{\mathrm{clo}(A)}\cK
\]
where $\mathrm{int}(A)$ and $\mathrm{clo}(A)$ are the interior and the closure
of $A$. Using the terminology of Srinivasa Varadhan\footnote{Srinivasa
  Varadhan (1940 -- ) is the (main) father of modern large deviations theory.},
${(\mu_N)}_{N\geq1}$ satisfies to a large deviations principle with speed $N$
and rate function $\cK$. The functional $\cK:\cM_1\to\dR\cup\{+\infty\}$ is
the Kullback-Leibler relative entropy with respect to $\mu_*$. By convention
$\cK(\mu):=+\infty$ if $\mu\not\ll\mu_*$. If
$d\mu_*(x)=f_*(x)\,dx=Z_\be^{-1}e^{-\be V}\,dx$ and $d\mu=f\,d\mu_*$ then
\[
\cK(\mu)=\be\PAR{\cA(f)-\cA(f_*)}.
\]
The Sanov theorem is a refinement of the strong law of large numbers, since by
the first Borel-Cantelli lemma, one obtains that with probability one,
$\lim_{N\to\infty}\mu_N=\mu_*=\arg\inf \cK$.

The large deviations rate function $\cK$ is convex and lower semicontinuous
with respect to the narrow topology, which is the topology of convergence with
respect to bounded and continuous test functions. This topology can be
metrized by the metric $d(\mu,\nu):=\sup_{h\in\cH}\int\!h\,d(\mu-\nu)$ where
$\cH:=\{h:\max(\NRM{h}_\infty,\LIP{h})\leq1\}$. Now for
$A=A_\veps=B(\mu,\veps):=\{\nu:d(\mu,\nu)\leq\veps\}$ we have
\[
\dP(\mu_N\in A)
=\mu_*^{\otimes
  N}\PAR{(x_1,\ldots,x_N)\in\dR^N:\sup_{h\in\cH}\PAR{\frac{1}{N}\sum_{i=1}^Nh(x_i)-\int\!h\,d\mu}\leq
\veps}.
\]
and thanks to the Sanov theorem, we obtain the ``volumetric'' formula
\[
\inf_{\veps>0}
\limsup_{N\to\infty}
\frac{1}{N}\log
\mu_*^{\otimes
  N}\PAR{(x_1,\ldots,x_N)\in\dR^N:\sup_{h\in\cH}\PAR{\frac{1}{N}\sum_{i=1}^Nh(x_i)-\int\!h\,d\mu}\leq
  \veps}
=-K(\mu).
\]

\subsection{Names}

The letter $\cS$ was chosen by Rudolf Clausius (1822 -- 1888) for entropy in
thermodynamics, possibly in honor of Sadi Carnot (1796 -- 1832). The term
\emph{entropy} was forged by Clausius in 1865 from the Greek «$\eta\
\tau\rho{o}\pi\eta$». The letter $H$ used by Boltzmann is the capital Greek
letter $\eta$. The letter $\cA$ used for the Helmholtz free energy comes from
the German word ``Arbeit'' for work.

\medskip

\begin{quote}``I propose to name the quantity $\cS$ the entropy of the system,
  after the Greek word $\eta\ \tau\rho{o}\pi\eta$ (en tropein), the
  transformation. I have deliberately chosen the word entropy to be as similar
  as possible to the word energy: the two quantities to be named by these
  words are so closely related in physical significance that a certain
  similarity in their names appears to be appropriate.''
  \begin{flushright}
    Rudolf Clausius, 1865
  \end{flushright}
\end{quote}

\subsection{H-Theorem}

Back to the motivations of Boltzmann, let us recall that the first principle
of Carnot-Clausius thermodynamics\footnote{According to Vladimir Igorevitch
  Arnold (1937 -- 2010), ``Every mathematician knows it is impossible to
  understand an elementary course in thermodynamics.''. Nevertheless, the
  reader may try \cite{fermi,MR1724307}, and \cite{MR1881344} for history.}
states that the internal energy of an isolated system is constant, while the
second principle states that there exists an extensive state variable called
the entropy that can never decrease for an isolated system. Boltzmann wanted
to derive the second principle from the idea (controversial, at that time)
that the matter is made with atoms. Let us consider an ideal isolated gas made
with particles (molecules) in a box with periodic boundary conditions (torus)
to keep things as simple as possible. There are of course too many particles
to write the equations of Newton for all of them. Newton is in a way beated by
Avogadro! The idea of Boltzmann was to propose a statistical approach (perhaps
inspired from the one of Euler in fluid mechanics, and from the work of
Maxwell, see \cite{MR1881344}): instead of keeping track of each particle, let
$(x,v)\mapsto f_t(x,v)$ be the probability density of the distribution of
position $x\in\dR^d$ and velocity $v\in\dR^d$ of particles at time $t$. Then
one can write an evolution equation for $t\mapsto f_t$, that takes into
account the physics of elastic collisions. It is a nonlinear partial
differential equation known as the Boltzmann equation:
\[
\pd_t f_t(x,v) =- v\pd_x f_t(x,v) + Q(f_t,f_t)(x,v).
\]
The first term in the right hand side is a linear transport term, while the
second term $Q(f_t,f_t)$ is quadratic in $f_t$, a double integral actually,
and captures the physics of elastic collisions by averaging over all possible
input and output velocities (note here a loss of microscopic information).
This equation admits conservation laws. Namely, for every time $t\geq0$, $f_t$
is a probability density and the energy of the system is constant (first
principle):
\[
f_t\geq0,\quad\pd_t\int\!f_t(x,v)\,dxdv=0,\quad
\pd_t\displaystyle\iint\!v^2\,f_t(x,v)\,dxdv=0.
\]
These constrains define a class of densities $\cC$ on $\dR^d\times\dR^d$ over
which the Boltzmann entropy $\cS$ achieves its (Gaussian in velocity and
uniform in position) maximum
\[
f_*=\arg\max_\cC \cS.
\]
The H-Theorem states that the entropy $\cS=-H$ is monotonic along the
Boltzmann equation:
\[
\pd_t \cS(f_t) \geq 0
\]
and more precisely,
\[
\cS(f_t)\underset{t\to\infty}{\nearrow}\cS(f_*)=\max_\cC \cS
\]
where $\cC$ is the class defined by the conservation law. A refined analysis
gives that 
\[
f_t\underset{t\to\infty}{\longrightarrow} f_*=\arg\max_\cC \cS.
\]
In the space-homogeneous simplified case, $f_t$ depends only on the velocity
variable, giving a Gaussian equilibrium for velocities by maximum entropy! In
kinetic theory of gases, it is customary to call ``Maxwellian law'' the
standard Gaussian law on velocities. We refer to \cite{MR3076094} for a
discussion on the concept of irreversibility and the Boltzmann H-Theorem.

The work of Boltzmann in statistical physics (nineteenth century) echoes the
works of Euler in fluid mechanics (eighteenth century), and of Newton in
dynamics (seventeenth century). Before the modern formalization of probability
theory and of partial differential equations with functional analysis,
Boltzmann, just like Euler, was able to forge a deep concept melting the two!
The Boltzmann H-Theorem had and has still a deep influence, with for instance
the works of Kac, Lanford, Cercignani, Sinai, Di Perna and Lions, Desvillettes
and Villani, Saint-Raymond, \ldots.

\medskip

\begin{quote}
  ``Although Boltzmann’s H-Theorem is 135 years old, present-day mathematics
  is unable to prove it rigorously and in satisfactory generality. The
  obstacle is the same as for many of the famous basic equations of
  mathematical physics: we don’t know whether solutions of the Boltzmann
  equations are smooth enough, except in certain particular cases (close-
  to-equilibrium theory, spatially homogeneous theory, close-to-vacuum
  theory). For the moment we have to live with this shortcoming.''
\end{quote}
\begin{flushright}
  Cédric Villani, 2008, excerpt from \cite{MR2509760}\\
  H-Theorem and beyond: Boltzmann’s entropy in today’s mathematics
\end{flushright}

\subsection{Keeping in mind the structure}
For our purposes, let us keep in mind this idea of evolution equation,
conservation law, monotonic functional, and equilibrium as optimum (of the
monotonic functional) under constraint (provided by the conservation law). It
will reappear!

\subsection{Markov processes and Helmholtz free energy} 

A Markov process can always be seen as a deterministic evolution equation of a
probability law. By analogy with the Boltzmann equation, let us consider a
Markov process $(X_t,t\in\dR_+)$ on $\dR^d$. Let us focus on structure and
relax the rigor to keep things simple. For any $t\geq0$ and continuous and
bounded test function $h$, for any $x\in\dR^d$,
\[
P_t(h)(x):=\dE(h(X_t)|X_0=x).
\]
Then $P_0(h)=h$, and, thanks to the Markov property, the one parameter family
$P=(P_t,t\in\dR_+)$ forms a semigroup of operators acting on bounded
continuous test functions, with $P_0=id$. Let us assume that the process
admits an invariant measure $\mu_*$, meaning that for every $t\geq0$ and $h$,
\[
\int\!P_t(h)\,d\mu_*=\int\!h\,d\mu_*.
\]
The semigroup is contractive in $L^p(\mu_*)$: $\NRM{P_t}_{p\to p}\leq1$ for
any $p\in[1,\infty]$ and $t\geq0$. The semigroup is Markov: $P_t(h)\geq0$ if
$h\geq0$, and $P_t(1)=1$. The infinitesimal generator of the semigroup is
\(
Lh=\pd_{t=0}P_t(h),
\)
for any $h$ in the domain of $L$ (we ignore these aspects for simplicity). In 
particular 
\[
\int\!Lh\,d\mu_*
=\pd_{t=0}\int\!P_t(h)\,d\mu_*
=\pd_{t=0}\int\!h\,d\mu_*
=0.
\]
When $L$ is a second order linear differential operator without constant term,
then we say that the process is a Markov diffusion. Important examples of such
Markov diffusions are given in Table \ref{tab:diffproc}. The backward and
forward Chapman-Kolmogorov equations are
\[
\pd_tP_t=LP_t=P_tL.
\]
Let us denote by $P_t^*$ the adjoint of $P_t$ in $L^2(\mu_*)$, and let us
define
\[
\mu_t:=\mathrm{Law}(X_t)\quad\text{and}\quad g_t:=\frac{d\mu_t}{d\mu_*}.
\]
Then 
$g_t=P_t^*(g_0)$. If we denote $L^*=\pd_{t=0}P_t^*$ then we obtain
the evolution equation
\[
\pd_tg_t=L^*g_t.
\]
The invariance of $\mu_*$ can be seen as a fixed point: if $g_0=1$ then
$g_t=P_t^*(1)=1$ for all $t\geq0$, and $L^*1=0$. One may prefer to express the
evolution of the density 
\[
f_t:=\frac{d\mu_t}{dx}=g_tf_*
\quad\text{where}\quad 
f_*:=\frac{d\mu_*}{dx}.
\]
We have then $\pd_tf_t=G^*f_t$ where $G^*h:=L^*(h/f_*)f_*$. 
The linear evolution equation
\[
\pd_tf_t=G^*f_t
\]
is in a sense the Markovian analogue of the Boltzmann equation (which is
nonlinear!).
\begin{table}
  \begin{tabular}{c||c|c|c}
    & Brownian motion 
    & Ornstein-Uhlenbeck process 
    & Overdamped Langevin process\\\hline
    S.D.E.\ 
    & $dX_t=\sqrt{2}dB_t$ 
    & $dX_t=\sqrt{2}dB_t-X_tdt$ 
    & $dX_t=\sqrt{2}dB_t-\nabla V(X_t)\,dt$\\
    $\mu_*$ 
    & $dx$ 
    & $\cN(0,I_d)$ 
    & $Z^{-1}e^{-V(x)}\,dx$\\
    $f_*$ 
    & $1$ 
    & $(2\pi)^{-d/2}e^{-\NRM{x}_2^2/2}$
    & $Z^{-1}e^{-V}$ \\
    $L f$ 
    & $\De f$ 
    & $\De f-x\cdot\nabla f$ 
    & $\De-\nabla V\cdot\nabla f$\\
    $G^* f$ 
    & $\De f$ 
    & $\De f+\mathrm{div}(xf)$
    & $\De+\mathrm{div}(f\nabla V)$\\
    $\mu_t$ 
    & $\mathrm{Law}(X_0+\sqrt{2t}G)$ 
    & $\mathrm{Law}(e^{-t}X_0+\sqrt{1-e^{-2t}}G)$
    & Not explicit in general\\ 
    $P_t$ 
    & Heat semigroup 
    & O.-U. semigroup 
    & General semigroup
  \end{tabular}
  \caption[Two fundamental Gaussian processes on $\dR^d$]{Two fundamental
    Gaussian processes on $\dR^d$, Brownian Motion and Ornstein-Uhlenbeck, as
    Gaussian special cases of Markov diffusion processes\protect\footnotemark.}
  \label{tab:diffproc}
\end{table}

\footnotetext{The $\sqrt{2}$ factor in the S.D.E.\ allows to avoid a factor
  $1/2$ in the infinitesimal generators.}

Let us focus on the case where $f_*$ is a density, meaning that $\mu_*$ is a
probability measure. This excludes Brownian motion, but allows the
Ornstein-Uhlenbeck process. By analogy with the Boltzmann equation, we have
the first two conservation laws $f_t\geq0$ and $\int\!f_t\,dx=1$, but the
average energy has no reason to be conserved for a Markov process. Indeed, the
following quantity has no reason to have a constant sign (one can check this
on the Ornstein-Uhlenbeck process!):
\[
\pd_t\int\!Vf_t\,dx
=\pd_tP_t(V)
=P_t(LV).
\]
Nevertheless, if we set $f_*=Z_\be^{-1}e^{-\be V}$ then there exists a
functional which is monotonic and which admits the invariant law $\mu_*$ as a
unique optimizer: the Helmholtz free energy defined by
\[
\cA(f):=\int\!V(x)\,f(x)\,dx-\frac{1}{\be}\cS(f).
\]
In order to compute $\pd_t\cA(f_t)$, we
first observe that for any test function $g$,
\[
\int\!L^*g\,d\mu_*
=0.
\]
Since $\mu_*$ is invariant, we have $P^*_t(1)=1$ for every $t\geq0$, and,
since $P_t^*(g)\geq0$ if $g\geq0$, it follows that the linear form $g\mapsto
P^*_t(g)(x)$ is a probability measure\footnote{Actually $P_t^*(g)(x)$ is the
  value at point $x$ of the density with respect to $\mu_*$ of the law of
  $X_t$ when $X_0\sim gd\mu_*$.}. Recall that
\[
\cA(f)-\cA(f_*)=\frac{1}{\be}\int\!\Phi(g)\,d\mu_*
\]
where $\Phi(u):=u\log(u)$. For every $0\leq s\leq t$, the Jensen inequality
and the invariance of $\mu_*$ give
\[
\cA(f_t)-\cA(f_*) %
=\int\!\Phi(P_{t-s}^*(f_s))\,d\mu_* %
\leq\int\!P_{t-s}^*(\Phi(f_s))\,d\mu_* %
=\int\!\Phi(f_s)\,d\mu_* %
=\cA(f_s)-\cA(f_*),
\]
which shows that the function $t\mapsto\cA(f_t)$ is monotonic. Alternatively,
the Jensen inequality gives also $\Phi(P^*_t(g))\leq P^*_t(\Phi(g))$ and the
derivative at $t=0$ gives $\Phi'(g)L^*g\leq L^*\Phi(g)$, which provides
\[
\int\!\Phi'(g)L^*g\,d\mu_*\leq0.
\]
Used with $g=g_t=f_t/f_*$, this gives
\[
\be\pd_t\cA(f_t)
=\pd_t\int\!\Phi(g_t)\,d\mu_*
=\int\!\Phi'(g_t)L^*g_t\,d\mu_*
\leq0.
\]
It follows that the Helmholtz free energy decreases along the Markov process:
\[
\pd_t\cA(f_t)\leq0.
\]
Of course we expect, possibly under more assumptions, that
$\cA(f_t)\searrow\cA(f_*)=\min\cA$ as $t\to\infty$. Let us assume for
simplicity that $\be=1$ and that the process is the Markov diffusion generated
by $L=\De -\nabla V\cdot\nabla$. In this case $\mu_*$ is symmetric for the
Markov process, and $L=L^*$, which makes most aspects simpler. By the chain
rule $L(\Phi(g))=\Phi'(g)Lg+\Phi''(g)\ABS{\nabla g}^2$, and thus, by
invariance,
\[
\pd_t\cA(f_t)
=\int\!\Phi'(g_t)Lg_t\,d\mu_*
=-\cF(g_t)\leq0
\quad\text{where}\quad
\cF(g):=
\int\!\Phi''(g)\ABS{\nabla g}^2\,d\mu_*
=\int\!\frac{\ABS{\nabla g}^2}{g}\,d\mu_*.
\]
The functional $\cF$ is known as the Fisher information\footnote{Named after
  Ronald Aylmer Fisher (1890 -- 1962), father of modern statistics among other
  things.}. The identity $\pd_t\cA(f_t)=-\cF(g_t)$ is known as the de Bruijn
identity. In the degenerate case $V\equiv0$, then $f_*\equiv1$ is no longer a
density, $\mu_*$ is the Lebesgue measure, the Markov process is Brownian
motion, the infinitesimal generator in the Laplacian $L=\De$, the semigroup
$(P_t,t\in\dR_+)$ is the heat semigroup, and we still have a de Bruijn
identity $\pd_t\cS(f_t)=\cF(f_t)$ where $\cS=-\cA$ since $V\equiv0$.

The quantitative version of the monotonicity of the free energy along the
Markov semigroup is related to Sobolev type functional inequalities. We refer
to \cite{MR1845806,BGL} for more details. For instance, for every constant
$\rho>0$, the following three properties are equivalent:
\begin{itemize}
\item Exponential decay of free energy: $\forall f_0\geq0$, $\forall t\geq0$,
  $\cA(f_t)-\cA(f_*)\leq e^{-2\rho t}(\cA(f_0)-\cA(f_*))$;
\item Logarithmic Sobolev inequality: $\forall f\geq0$,
  $2\rho(\cA(f)-\cA(f_*))\leq \cF(f/f_*)$;
\item Hypercontractivity: $\forall t\geq0$, $\NRM{P_t}_{q(0)\to q(t)}\leq1$
  where $q(t):=1+e^{2\rho t}$.
\end{itemize}
The equivalence between the first two properties follows by taking the
derivative over $t$ and by using the Grönwall lemma. The term ``Logarithmic
Sobolev inequality'' is due to Leonard Gross, who showed in \cite{MR420249}
the equivalence with hypercontractivity, via the basic fact that for any
$g\geq0$,
\[
\pd_{p=1}\NRM{g}_p^p
=\pd_{p=1}\int\!e^{p\log(g)}\,d\mu_*
=\int\!g\log(g)\,d\mu_*
=\cA(gf_*)-\cA(f_*).
\]
The concept of hypercontractivity of semigroups goes back at least to Edward
Nelson \cite{MR0210416}.

One may ask if $t\mapsto\cF(f_t)$ is in turn monotonic. The answer involves a
notion of curvature. Namely, using the diffusion property via the chain rule
and reversibility, we get, after some algebra,
\[
\pd_t^2\cA(f_t)
=-\pd_t\cF(g_t)
=2\!\!\int\!g_t\Ga_{\!\!2}(\log(g_t))\,d\mu_*
\]
where $\Ga_{\!\!2}$ is the Bakry-Émery ``gamma-two'' functional quadratic form
given by\footnote{The terminology comes from
  $\Ga_{\!\!2}(f):=\Ga_{\!\!2}(f,f):=\frac{1}{2}(L(\Ga(f,f))-2\Ga(f,Lf))$
  where $\Ga$ is the ``carré du champ'' functional quadratic form
  $(f,g)\mapsto\Ga(f,g)$ defined by
  $\Ga(f)=\Ga(f,f):=\frac{1}{2}(L(f^2)-2fLf)$. Here $\Ga(f)=\ABS{\nabla
    f}^2$.}
\[
\Ga_{\!\!2}(f)
=\NRM{\nabla^2 f}_{\mathrm{HS}}^2
+\nabla f\cdot (\nabla^2 V)\nabla f.
\]
See \cite{MR1845806,BGL} for the details. This comes from the Bochner
commutation formula
\[
\nabla L=L\nabla-(\nabla^2V)\nabla.
\]
If $\Ga_{\!\!2}\geq0$ then, along the semigroup, the Fisher information is
non-increasing, and the Helmholtz free energy is convex (we already know that
$\pd_t\cA(f_t)\leq0$ and $\min\cA=\cA(f_*)$):
\[
\pd_t\cF(g_t)\leq0\quad\text{and}\quad\pd^2_t\cA(f_t)\geq0.
\]
This holds for instance if $\nabla^2 V\geq0$ as quadratic forms\footnote{This
  means $y\cdot(\nabla^2 V(x))y\geq0$ for all $x,y\in\dR^d$, or equivalently
  $V$ is convex, or equivalently $\mu_*$ is log-concave.}. This is the case
for the Ornstein-Uhlenbeck example, for which $\nabla^2 V=I_d$. Moreover, if
there exists a constant $\rho>0$ such that for all $f$,
\[
\Ga_{\!\!2}(f)\geq\rho\Ga(f) \quad\text{where}\quad\Ga(f)=\ABS{\nabla f}^2
\]
then, for any $t\geq0$, $\pd_t\cF(g_t)\leq-2\rho\cF(g_t)$, and the Grönwall
lemma gives the exponential decay of the Fisher information along the
semigroup ($\rho=1$ in the Ornstein-Uhlenbeck example):
\[
\cF(g_t)\leq e^{-2\rho t}\cF(g_0).
\]
This gives also the exponential decay of the Helmholtz free energy $\cA$ along
the semigroup with rate $2\rho$, in other words, a Logarithmic Sobolev
inequality with constant $2\rho$: for any $f_0$,
\[
\cA(f_0)-\cA(f_*)
=-\int_0^\infty\!\pd_t\cA(f_t)\,dt
=\int_0^\infty\!\cF(g_t)\,dt
\leq \cF(g_0)\int_0^\infty\!e^{-2\rho t}\,dt
=\frac{\cF(f_0/f_*)}{2\rho}.
\]
We used here the Markov semigroup in order to interpolate between $f_0$ and
$f_*$. This interpolation technique was extensively developed by Bakry and
Ledoux in order to obtain functional inequalities, see \cite{BGL} and
references therein. A modest personal contribution to this topic is
\cite{MR2081075}. The best possible constant $\rho$ -- which is the largest --
in the inequality $\Ga_{\!\!2}\geq\rho\Ga$ is called the Bakry-Émery curvature
of the Markov semigroup. The story remains essentially the same for Markov
processes on Riemannian manifolds. More precisely, in the context of Riemannian
geometry, the Ricci curvature tensor contributes additively to the
$\Ga_{\!\!2}$, see \cite{BGL}. Relatively recent works on this topic include
extensions by Lott and Villani \cite{MR2480619}, von Renesse and Sturm
\cite{MR2142879}, Ollivier \cite{MR2484937}, among others. The approach can be
adapted to non-elliptic hypoelliptic evolution equations, see for instance
Baudoin \cite{baudoin}. The Boltzmannian idea of monotonicity along an
evolution equation is also used in the work of Grigori Perelman on the
Poincaré-Thurston conjecture, and in this case, the evolution equation, known
as the Ricci flow of Hamilton, concerns the Ricci tensor itself, see
\cite{perelman,xdli}. The exponential decay of the Boltzmann functional
$H=-\cS$ along the (nonlinear!) Boltzmann equation was conjectured by
Cercignani and studied rigorously by Villani, see for instance
\cite{MR2765747}.

Many aspects remain valid for discrete time/space Markov chains, up to the
lack of chain rule if the space is discrete. For instance, if an irreducible
Markov chain in discrete time and finite state space $E$ has transition matrix
$P$ and invariant law $\mu_*$, then ${(P^n)}_{n\geq0}$ is the discrete time
Markov semigroup, $L:=P-I$ is the Markov generator, and one can show that for
every initial law $\mu$, the discrete Helmholtz free energy is monotonic along
the evolution equation, namely
\[
\sum_{x\in E}\Phi\PAR{\frac{\mu P^n(x)}{\mu_*(x)}}\mu_*(x)
\underset{n\to\infty}{\searrow}
0,
\]
where $\mu P^n(x)=\sum_{z\in E}\mu(z)P^n(z,x)$ and still $\Phi(u):=u\log(u)$.
The details are in Thomas Liggett's book \cite[prop.~4.2]{MR2108619}.
According to Liggett \cite[p.~120]{MR2108619} this observation goes back at
least to Mark Kac \cite[p.~98]{MR0102849}. Quoting Persi Diaconis, ``\emph{the
  idea is that the maximum entropy Markov transition matrix with a given
  invariant law is clearly the matrix with all rows equal to this stationary
  law, taking a step in the chain increases entropy and keeps the stationary
  law the same.}''.

Discrete versions of the logarithmic Sobolev inequality allow to refine the
quantitative analysis of the convergence to the equilibrium of finite state
space Markov chains. We refer for these aspects to the work of Diaconis and
Saloff-Coste \cite{MR1410112,MR1490046}, the work of Laurent Miclo
\cite{MR1478724}, and the book \cite{MT}. The relative entropy allows to
control the total variation distance: the so-called Pinsker or
Csiszár-Kullback inequality states that for any probability measures $\mu$ and
$\nu$ on $E$ with $\mu>0$,
\[
\sum_{x\in E}\ABS{\mu(x)-\nu(x)}
\leq \sqrt{2\sum_{x\in E}\Phi\PAR{\frac{\nu(x)}{\mu(x)}}\mu(x)}.
\]
The analysis of the convexity of the free energy along the semigroup of at
most countable state space Markov chains was considered in \cite{CDPP} and
references therein. More precisely, let ${(X_t)}_{t\in\dR_+}$ be a continuous
time Markov chain with at most countable state space $E$. Let us assume that
it is irreducible, positive recurrent, and aperiodic, with unique invariant
probability measure $\mu_*$, and with infinitesimal generator $L:E\times
E\to\dR$. We have, for every $x,y\in E$,
\[
L(x,y)=\pd_{t=0}\dP(X_t=y|X_0=x).
\]
We see $L$ as matrix with non-negative off-diagonal elements and zero-sum
rows: $L(x,y)\geq0$ and $L(x,x)=-\sum_{y\neq x}L(x,y)$ for every $x,y\in E$.
The invariance reads $0=\sum_{x\in E}\mu_*(x)L(x,y)$ for every $y\in E$. The
operator $L$ acts on functions as $(Lf)(x)=\sum_{y\in E}L(x,y)f(y)$ for every
$x\in E$. Since $\mu_*(x)>0$ for every $x\in E$, the free energy at unit
temperature corresponds to the energy $V(x)=-\log(\mu_*(x))$, for which we
have of course $\cA(\mu_*)=0$. For any probability measure $\mu$ on $E$,
\[
\cA(\mu)-\cA(\mu_*)
=\cA(\mu)
=\sum_{x\in E}\Phi\PAR{\frac{\mu(x)}{\mu_*(x)}}\mu_*(x).
\]
One can see $x\mapsto\mu(x)$ as a density with respect to the counting measure
on $E$. For any time $t\in\dR_+$, if $\mu_t(x):=\dP(X_t=x)$ then
$g_t(x):=\mu_t(x)/\mu_*(x)$ and $\pd_t g_t=L^*g_t$ where $L^*$ is the adjoint
of $L$ in $\ell^2(\mu_*)$ which is given by
$L^*(x,y)=L(y,x)\mu_*(y)/\mu_*(x)$. Some algebra reveals that
\[
\pd_t\cA(\mu_t)=\sum_{x\in E}\SBRA{\Phi'(g_t)L^*g_t}(x)\mu_*(x).
\]
The right hand side is up to a sign the discrete analogue of the Fisher
information. By reusing the convexity argument used before for diffusions, we
get that $\pd_t\cA(\mu_t)\leq0$. Moreover, we get also
\[
\pd_t^2\cA(\mu_t)
=\sum_{x\in E}\SBRA{g_tLL\log(g_t)+\frac{(L^*g_t)^2}{g_t}}(x)\mu_*(x).
\]
The right hand side is a discrete analogue of the $\Ga_{\!\!2}$-based formula
obtained for diffusions. It can be nicely rewritten when $\mu_*$ is
reversible. The lack of chain rule in discrete spaces explains the presence of
two distinct terms in the right hand side. We refer to \cite{CDPP,venkat} for
discussions of examples including birth-death processes. Our modest
contribution to this subject can be found in \cite{MR2247924,MR3129037}. We
refer to \cite{MR2484937} for some complementary geometric aspects. An
analogue on $E=\dN$ of the Ornstein-Uhlenbeck process is given by the
so-called M/M/$\infty$ queue for which
$(Lf)(x)=\la(f(x+1)-f(x))+x\mu((f(x-1)-f(x))$ and
$\mu_*=\mathrm{Poisson}(\la/\mu)$.


\section{Claude Shannon and the central limit theorem}

The Boltzmann entropy plays also a fundamental role in communication theory,
funded in the 1940's by Claude Elwood Shannon (1916--2001), where it is known
as ``Shannon entropy''. It has a deep interpretation in terms of uncertainty
and information in relation with coding theory \cite{MR2239987}. For example
the discrete Boltzmann entropy $\cS(p)$ computed with a logarithm in base $2$ is
the average number of bits per symbol needed to encode a random text with
frequencies of symbols given by the law $p$. This plays an essential role in
lossless coding, and the Huffman algorithm for constructing Shannon entropic
codes is probably one of the most used basic algorithm (data compression is
everywhere). Another example concerns the continuous Boltzmann entropy which
enters the computation of the capacity of continuous telecommunication
channels (e.g.\ DSL lines).

\begin{quote}
  ``My greatest concern was what to call it. I thought of calling it
  ‘information’, but the word was overly used, so I decided to call it
  ‘uncertainty’. When I discussed it with John von Neumann, he had a better
  idea. Von Neumann told me, ‘You should call it entropy, for two reasons. In
  the first place your uncertainty function has been used in statistical
  mechanics under that name, so it already has a name. In the second place,
  and more important, nobody knows what entropy really is, so in a debate you
  will always have the advantage.''
\end{quote}
\begin{flushright}
  Claude E. Shannon, 1961\\
  Conversation with Myron Tribus, reported in \cite{Tribus:1971:EI}
\end{flushright}

For our purposes, let us focus on the link between the Boltzmann entropy and
the central limit theorem, a link suggested by Shannon when he forged
information theory in the 1940's.

\subsection{The CLT as an evolution equation}

The Central Limit theorem (CLT) states that if $X_1,X_2,\ldots$ are i.i.d.\
real random variables with mean $\dE(X_i)=0$ and variance $\dE(X_i^2)=1$, then
\[
S_n:= \frac{X_1+\cdots+X_n}{\sqrt{n}}
\underset{n\to\infty}{\overset{d}{\longrightarrow}}
\frac{e^{-\frac{1}{2}x^2}}{\sqrt{2\pi}}dx
\]
where the convergence to the Gaussian law holds in distribution (weak sense).

\subsection{Conservation law}

The first two moments are conserved along CLT: for all $n\geq1$,
\[
\dE(S_n)=0\quad\text{and}\quad\dE(S_n^2)=1.
\]
By analogy with the H-Theorem, the CLT concerns an evolution equation of the
law of $S_n$ along the discrete time $n$. When the sequence $X_1,X_2,\ldots$
is say bounded in $L^\infty$ then the convergence in the CLT holds in the
sense of moments, and in particular, the first two moments are constant while
the remaining moments of order $>2$ become universal at the limit. In other
words, the first two moments is the sole information retained by the CLT from
the initial data, via the conservation law. The limiting distribution in the
CLT is the Gaussian law, which is the maximum of Boltzmann entropy under
second moment constraint. If we denote by $f^{*n}$ the $n$-th convolution
power of the density $f$ of the $X_i$'s then the CLT writes
$\mathrm{dil}_{n^{-1/2}}(f^{*n})\to f_*$ where $f_*$ is the standard Gaussian
and where $\mathrm{dil}_\al(h):=\al^{-1}h(\al^{-1}\cdot)$ is the density of
the random variable $\al Z$ when $Z$ has density $h$.

\subsection{Analogy with H-Theorem}

Shannon observed \cite{MR0032134} that the entropy $\cS$ is monotonic along
the CLT when $n$ is a power of $2$, in other words $\cS(S_{2^{m+1}})\geq
\cS(S_{2^{m}})$ for every integer $m\geq0$, which follows from (a rigorous
proof is due to Stam \cite{stam})
\[
\cS\PAR{\frac{X_1+X_2}{\sqrt{2}}}=\cS(S_2)\geq \cS(S_1)=\cS(X_1).
\]
By analogy with the Boltzmann H-Theorem, a conjecture attributed to Shannon
(see also \cite{MR0506364}) says that the Boltzmann entropy $\cS$ is monotonic
along the CLT for any $n$, more precisely
\[
\cS(X_1)=\cS(S_1)\leq\cdots\leq \cS(S_n)\leq \cS(S_{n+1})
\leq\cdots\underset{n\to\infty}{\nearrow}\cS(G).
\]
The idea of proving the CLT using the Boltzmann entropy is very old and goes
back at least to Linnik \cite{MR0124081} in the 1950's, who, by the way, uses
the term ``Shannon entropy''. But proving convergence differs from proving
monotonicity, even if these two aspects are obviously lin(ni)ked
\cite{MR2128238}. The approach of Linnik was further developed by Rényi,
Csiszár, and many others, and we refer to the book of Johnson \cite{MR2109042}
for an account on this subject. The first known proof of the Shannon
monotonicity conjecture is relatively recent and was published in 2004 by
Artstein, Ball, Barthe, Naor \cite{MR2083473}. The idea is to pull back the
problem to the monotonicity of the Fisher information. Recall that the Fisher
information of a random variable $S$ with density $g$ is given by
\[
\cF(S):=\int\!\frac{\ABS{\nabla g}^2}{g}\,dx.
\]
It appears when one takes the derivative of $\cS$ along an additive Gaussian
perturbation. Namely, the de Bruijn formula states that if $X,G$ are random
variables with $G$ standard Gaussian then
\[
\pd_t\cS(X+\sqrt{t}G)=\frac{1}{2}\cF(X+\sqrt{t}G).
\]
Indeed, if $f$ is the density of $X$ then the density $P_t(f)$ of
$X+\sqrt{t}G$ is given by the heat kernel
\[
P_t(f)(x)=(f*\mathrm{dil}_{\sqrt{t}} f_*)(x)
=\int\!f(y)\frac{e^{-\frac{1}{2t}(y-x)^2}}{\sqrt{2\pi t}}\,dy,
\]
which satisfies to $\pd_tP_t(f)(x)=\frac{1}{2}\Delta_x P_t(f)(x)$, and which
gives, by integration by parts,
\[
\pd_t\cS(X+\sqrt{t}G)
=-\frac{1}{2}\int\!(1+\log P_t f)\De P_t\,dx
=\frac{1}{2}\cF(P_tf).
\]
In the same spirit, we have the following integral representation, taken from
\cite{MR2346565},
\[
\cS(G)-\cS(S)%
=\int_0^{\infty}\!
\PAR{\cF(\sqrt{e^{-2t}}S+\sqrt{1-e^{-2t}}G)-1}\,dt.
\]  
This allows to deduce the monotonicity of $\cS$ along the CLT from the one of
$F$ along the CLT,  
\[
\cF(S_1)\geq \cF(S_2)\geq\cdots\geq \cF(S_n)\geq
\cF(S_{n+1})\geq\cdots\searrow \cF(G),
\]
which is a tractable task \cite{MR2083473}. The de Bruijn identity involves
Brownian motion started from the random initial condition $X$, and the
integral representation above involve the Ornstein-Uhlenbeck process started
from the random initial condition $X$. The fact that the Fisher information
$\cF$ is non-increasing along a Markov semigroup means that the entropy $\cS$ is
concave along the Markov semigroup, a property which can be traced back to
Stam (see \cite[Chapter 10]{MR1845806}). As we have already mentioned before,
the quantitative version of such a concavity is related to Sobolev type
functional inequalities and to the notion of Bakry-Émery curvature of Markov
semigroups \cite{BGL}, a concept linked with the Ricci curvature in Riemannian
geometry. Recent works on this topic include extensions by Villani, Sturm,
Ollivier, among others.

\section{Dan-Virgil Voiculescu and the free central limit theorem}

Free probability theory was forged in the 1980's by Dan-Virgil Voiculescu
(1946--), while working on isomorphism problems in von Neumann operator
algebras of free groups. Voiculescu discovered later in the 1990's that free
probability is the algebraic structure that appears naturally in the
asymptotic global spectral analysis of random matrix models as the dimension
tends to infinity. Free probability theory comes among other things with
algebraic analogues of the CLT and the Boltzmann entropy, see
\cite{MR1217253,MR1887698,MR2032363,AGZ}. The term ``free'' in ``free
probability theory'' and in ``free entropy'' comes from the free group (see
below), and has no relation with the term ``free'' in the Helmholtz free
energy which comes from thermodynamics (available work obtainable at constant
temperature). By analogy, the ``free free energy'' at unit temperature might
be $\cA_*(a)=\tau(V(a))- \chi(a)$ where $\chi$ is the Voiculescu free entropy.
We will see in the last sections that such a functional appears as the rate
function of a large deviations principle for the empirical spectral
distribution of random matrix models! This is not surprising since the
Helmholtz free energy, which is nothing else but a Kullback-Leibler relative
entropy, is the rate function of the large deviations principle of Sanov,
which concerns the empirical measure version of the law of large numbers.

\subsection{Algebraic probability space}

Let $\cA$ be an algebra over $\dC$, with unity $id$, equipped with an
involution $a\mapsto a^*$ and a normalized linear form $\tau:\cA\to\dC$ such
that $\tau(ab)=\tau(ba)$, $\tau(id)=1$, and $\tau(aa^*)\geq0$. A basic non
commutative example is given by the algebra of square complex matrices:
$\cA=\cM_n(\dC)$, $id=I_n$, $a^*=\bar{a}^\top$, $\tau(a)=\frac{1}{n}\TR(a)$,
for which $\tau$ appears as an expectation with respect to the empirical
spectral distribution: denoting $\la_1(a),\ldots,\la_n(a)\in\dC$ the
eigenvalues of $a$, we have
\[
\tau(a)=\frac{1}{n}\sum_{k=1}^n\la_k(a)=\int\!x\,d\mu_a(x)
\quad\text{where}\quad
\mu_a:=\frac{1}{n}\sum_{k=1}^n\de_{\la_k(a)}.
\]
If $a=a^*$ (we say that $a$ is real or Hermitian) then the probability measure
$\mu_a$ is supported in $\dR$ and is fully characterized by the collection of
moments $\tau(a^m)$, $m\geq0$, which can be seen as a sort of algebraic
distribution of $a$. Beyond this example, by analogy with classical
probability theory, we may see the elements of $\cA$ as algebraic analogues of
bounded random variables, and $\tau$ as an algebraic analogue of an
expectation, and $\tau(a^m)$ as the analogue of the $m$-th moment of $a$. We
say that $a\in\cA$ has mean $\tau(a)$ and variance
$\tau((a-\tau(a))^2)=\tau(a^2)-\tau(a)^2$, and that $a$ is centered when
$\tau(a)=0$. The $*$-law of $a$ is the collection of mixed moments of $a$ and
$a^*$ called the $*$-moments:
\[
\tau(b_1\cdots b_m)\text{ where } b_1,\ldots,b_m\in\{a,a^*\}\text{ and } m\geq1.
\]		
In contrast with classical probability, the product of algebraic variables may
be non commutative. When $a\in\cA$ is real $a^*=a$, then the $*$-law of $a$
boils down to the moments: $\tau(a^m)$, $m\geq0$. In classical probability
theory, the law of a real bounded random variable $X$ is characterized by its
moments $\dE(X^m)$, $m\geq0$, thanks to the (Stone-)Weierstrass theorem. When
the bounded variable is not real and takes its values in $\dC$ then we need
the mixed moments $\dE(X^m\bar{X}^n)$, $m,n\geq0$.

One can connect $*$-law and spectrum even for non real elements. Namely, if
$a\in\cM_n(\dC)$ and if $\mu_a:=\frac{1}{n}\sum_{k=1}^n\de_{\la_k(a)}$ is its
empirical spectral distribution in $\dC$, then, for any
$z\not\in\{\la_1(a),\ldots,\la_n(a)\}$,
\begin{align*}
  \frac{1}{2}\tau(\log((a-zid)(a-zid)^*))
  &=\frac{1}{n}\log\ABS{\det(a-zI_n)} \\
  &=\int\!\log\ABS{z-\la}\,d\mu_a(\la) \\
  &=(\log\ABS{\cdot}*\mu_a)(z) \\
  &=:-U_{\mu_a}(z).
\end{align*}
The quantity $U_{\mu_a}(z)$ is exactly the logarithmic potential at point
$z\in\dC$ of the probability measure $\mu_a$. Since
$-\frac{1}{2\pi}\log\ABS{\cdot}$ is the so-called fundamental solution of the
Laplace equation in dimension $2$, it follows that in the sense of Schwartz
distributions,
\[
\mu_a=\frac{1}{2\pi}\De U_{\mu_a}.
\]
Following Brown, beyond the matrix case $\cA=\cM_n(\dC)$, this suggest to
define the spectral measure of an element $a\in\cA$ in a abstract algebra
$\cA$ as being the probability measure $\mu_a$ on $\dC$ given by
\[
\mu_a:=-\frac{1}{\pi}\De\tau(\log((a-zid)(a-zid)^*))
\]
where here again $\De=\pd\OL{\pd}$ is the two-dimensional Laplacian acting on
$z$, as soon as we know how to define the operator $\log((a-zid)(a-zid)^*)$
for every $z$ such that $a-zid$ is invertible. This makes sense for instance
if $\cA$ is a von Neumann algebra of bounded operators on a Hilbert space,
since one may define $\log(b)$ if $b^*=b$ by using functional calculus. The
moral of the story is that the $*$-law of $a$ determines the $*$-law of the
Hermitian element $(a-zid)(a-zid)^*$ for every $z\in\dC$, which in turn
determines the Brown spectral measure $\mu_a$. This strategy is known as
Hermitization.

The so-called Gelfand-Naimark-Segal (GNS, see \cite{AGZ}) construction shows
that any algebraic probability space can be realized as a subset of the
algebra of bounded operators on a Hilbert space. Using then the spectral
theorem, this shows that any compactly supported probability measure is the
$*$-law of some algebraic random variable.

\subsection{Freeness}

The notion of freeness is an algebraic analogue of independence. In classical
probability theory, a collection of $\si$-field are independent if the product
of bounded random variables is centered as soon as the factors are centered
and measurable with respect to different $\si$-fields. We say that
$\cB\subset\cA$ is a sub-algebra of $\cA$ when it is stable by the algebra
operations, by the involution, and contains the unity $id$. By analogy with
classical probability, we say that the collection ${(\cA_i)}_{i\in I}$ of
sub-algebras of $\cA$ are \emph{free} when for any integer $m\geq1$,
$i_1,\ldots,i_m\in I$, and $a_1\in\cA_{i_1},\ldots,a_m\in\cA_{i_m}$, we have
\[
\tau(a_1\cdots a_m)=0
\]
as soon as $\tau(a_1)=\cdots=\tau(a_m)=0$ and $i_1\neq \cdots\neq i_m$ (only
consecutive indices are required to be different). We say that ${(a_i)}_{i\in
  I}\subset\cA$ are free when the sub-algebras that they generate are free. If
for instance $a,b\in\cA$ are free and centered, then $\tau(ab)=0$, and
$\tau(abab)=0$. Note that in classical probability, the analogue of this last
expression will never be zero if $a$ and $b$ are not zero, due to the
commutation relation $abab=a^2b^2$.

Can we find examples of free matrices in $\cM_n(\dC)$? Actually, this will not
give exciting answers. Freeness is more suited for infinite dimensional
operators. It turns out that the definition and the name of freeness come from
a fundamental infinite dimensional example constructed from the free group.
More precisely, let $F_n$ be the free group\footnote{``Free'' because it is
  the free product of $n$ copies of $\mathbb{Z}$, without additional
  relation.} with $2n$ generators (letters and anti-letters)
$g_1^{\pm1},\ldots,g_n^{\pm1}$ with $n\geq2$ (for $n=1$, $F_n=\dZ$ is
commutative). Let $\varn$ be the neutral element of $F_n$ (empty string). Let
$A_n$ be the associated free algebra identified with a sub-algebra of
$\ell_\dC^2(F_n)$. Each element of $A_n$ can be seen as a finitely supported
complex measure of the form $\sum_{w\in F_n}c_w\de_w$. The collection
${(\de_w)}_{w\in F_n}$ is the canonical basis:
$\DOT{\de_w,\de_{w'}}=\IND_{w=w'}$. The product on $A_n$ is the convolution of
measures on $F_n$:
\[
\PAR{\sum_{w\in F_n}c_w\de_w}\PAR{\sum_{w\in F_n}c'_w\de_w} %
=\sum_{w\in F_n}\PAR{\sum_{v\in F_n}c_vc'_{v^{-1}w}}\de_w %
=\sum_{w\in F_n}(c*c')_w\de_w.
\]
Now, let $\cA$ be the algebra over $\dC$ of linear operators from
$\ell_\dC^2(F_n)$ to itself. The product in $\cA$ is the composition of
operators. The involution $*$ in $\cA$ is the transposition-conjugacy of
operators. The identity operator is denoted $id$. We consider the linear form
$\tau:\cA\to\dC$ defined for every $a\in\cA$ by
\[
\tau(a)=\DOT{a\de_\varn,\de_\varn}.
\]
For every $w\in F_n$, let $u_w\in\cA$ be the left translation operator defined
by $u_w(\de_v)=\de_{wv}$. Then 
\[
u_w u_{w'}=u_{ww'},\quad (u_w)^*=u_{w^{-1}},\quad u_\varn=id,
\]
and therefore $u_w u_w^*=u_w^* u_w=id$ (we say that $u_w$ is unitary). We have
$\tau(u_w)=\IND_{w=\varn}$. Let us show that the sub-algebras
$\cA_1,\ldots,\cA_n$ generated by $u_{g_1},\ldots,u_{g_n}$ are free. Each $\cA_i$
consists in linear combinations of $u_{g_i^r}$ with $r\in\dZ$, and centering
forces $r=0$. For every $w_1,\ldots,w_m\in F_n$, we have
\[
\tau(u_{w_1}\cdots u_{w_m})%
=\tau(u_{w_1\cdots w_m})%
=\IND_{w_1\cdots w_m=\varn}.
\]
Let us consider the case $w_j=g_{i_j}^{r_j}\in\cA_{i_j}$ with $r_j\in\dZ$, for
every $1\leq j\leq m$. Either $r_j=0$ and $w_j=\varn$ or $r_j\neq0$ and
$\tau(w_j)=0$. Let us assume that $r_1\neq0,\ldots,r_n\neq0$, which implies
$\tau(u_{w_1})=\cdots=\tau(u_{w_n})=0$.
Now since $G_n$ is a tree, it does not have cycles, and thus if we follow a
path starting from the root $\varn$ then we cannot go back to the root if we
never go back locally along the path (this is due to the absence of cycles).
Consequently, if additionally $i_1\neq\cdots\neq i_m$ (i.e.\ two consecutive
terms are different), then we have necessarily $w_1\cdots w_m\neq\varn$, and
therefore $\tau(u_{w_1}\cdots u_{w_m})=0$. From this observation one can
conclude that $\cA_1,\ldots,\cA_n$ are free.


Beyond the example of the free group: if $G_1,\ldots,G_n$ are groups and $G$
their free product, then the algebras generated by $\{u_g:g\in
G_1\},\ldots,\{u_g:g\in G_n\}$ are always free in the one generated by
$\{u_g:g\in G\}$.

\begin{table}
  \begin{tabular}{c|c}
    Classical probability 
    & Free probability\\\hline
    Bounded r.v.\ $X$ on $\dC$ 
    & Algebra element $a\in\cA$\\
    $\dE(X^m\bar{X}^{n})$ 
    & $\tau(b_1\cdots b_m)$, $b\in\{a,a^*\}^m$ \\
    Law = Moments
    & Law = $*$-moments\\
    $X$ is real 
    & $a=a^*$ \\
    Independence 
    & Freeness\\
    Classical convolution $*$
    & Free convolution $\boxplus$\\
    Gaussian law (with CLT)
    & Semicircle law (with CLT)\\
    Boltzmann entropy $\cS$
    & Voiculescu entropy $\chi$
  \end{tabular}
  \caption{Conceptual dictionary between classical probability and free
    probability. The first is 
    commutative while the second is typically non commutative. Free
    probability is the algebraic structure that emerges from the asymptotic
    analysis, over the dimension, of the empirical spectral distribution of
    unitary invariant random matrices. The term free comes from the algebra of
    linear operators over the free algebra of the free group, an example in
    which the concept of freeness emerges naturally, in relation with the
    symmetric random walk on the infinite regular tree of even degree $\geq 4$
    (which is the Cayley graph of a non-commutative free group).} 
\end{table}

\subsection{Law of free couples and free convolution}

In classical probability theory, the law of a couple of independent random
variables is fully characterized by the couple of laws of the variables. In
free probability theory, the $*$-law of the couple $(a,b)$ is the collection
of mixed moments in $a,a^*,b,b^*$. If $a,b$ are free, then one can compute the
$*$-law of the couple $(a,b)$ by using the $*$-law of $a$ and $b$, thanks to
the centering trick. For instance, in order to compute $\tau(ab)$, we may
write using freeness $0=\tau((a-\tau(a))(b-\tau(b)))=\tau(ab)-\tau(a)\tau(b)$
to get $\tau(ab)=\tau(a)\tau(b)$. As a consequence, one can show that the
$*$-law of a couple of free algebraic variables is fully characterized by the
couple of $*$-laws of the variables. This works for arbitrary vectors of
algebraic variables.

In classical probability theory, the law of the sum of two independent random
variables is given by the convolution of the law of the variables. In free
probability theory, the $*$-law of the sum $a+b$ of two free algebraic
variables $a,b\in\cA$ is given by the so-called free convolution
$\mathrm{dist}(a)\boxplus\mathrm{dist}(b)$ of the $*$-law of $a$ and $b$,
which can be defined using the $*$-law of the couple $(a,b)$. Following
\cite{MR2032363}, given an at most countable family of compactly supported
probability measures on $\dR$, one can always construct an algebraic
probability space containing free algebraic variables admitting these
probability measures as their $*$-distributions. The free convolution
$\boxplus$ of probability measures is associative but is not distributive with
respect to convex combinations (beware of mixtures!).

We have so far at hand an algebraic framework, called free probability theory,
in which the concepts of algebra elements, trace, $*$-law, freeness, and free
convolution are the analogue of the concepts of bounded random variables,
expectation, law, independence, and convolution of classical probability
theory. Do we have a CLT, and an analogue of the Gaussian? The answer is
positive.

\subsection{Free CLT and semicircle law}
It is natural to define the convergence in $*$-law, denoted
$\overset{*}{\to}$, as being the convergence of all $*$-moments. The
Voiculescu free CLT states that if $a_1,a_2,\ldots\in\cA$ are free, real
$a_i=a_i^*$, with same $*$-law, zero mean $\tau(a_i)=0$, unit
variance $\tau(a_i^2)=1$, then
\[
s_n:=\frac{a_1+\cdots+a_n}{\sqrt{n}}
\overset{*}{\underset{n\to\infty}{\longrightarrow}}
\frac{\sqrt{4-x^2}\mathbf{1}_{[-2,2]}}{2\pi}dx.
\]
The limiting $*$-law is given by the moments of the \emph{semicircle
  law}\footnote{Also known as the Wigner distribution (random matrix theory)
  or the Sato-Tate distribution (number theory).} on $[-2,2]$, which are $0$
for odd moments and the Catalan numbers ${(C_m)}_{m\geq0}$ for even moments:
for every $m\geq0$,
\[
\int_{-2}^2\!x^{2m+1}\frac{\sqrt{4-x^2}}{2\pi}dx=0
\quad\text{and}\quad
\int_{-2}^2\!x^{2m}\frac{\sqrt{4-x^2}}{2\pi}dx
=C_m:=\frac{1}{1+m}\binom{2m}{m}.
\]
An algebraic variable $b\in\cA$ has semicircle $*$-law when it is real $b=b^*$
and $\tau(b^{2m+1})=0$ and $\tau(b^{2m})=C_m$ for every $m\geq0$. The proof of
the free CLT consists in computing the moments of $s_n$ using freeness. This
reveals three type of terms: terms which are zero at fixed $n$ thanks to
freeness and centering, terms having zero contribution asymptotically as
$n\to\infty$, and terms which survive at the limit, and which involve only the
second moment of the $a_i$'s. See \cite{MR2316893}\footnote{Beyond freeness,
  there exists other notions of algebraic independence, allowing to compute
  moments, and leading to a CLT with other limiting distributions, including
  the Bernoulli distribution and the arcsine distribution.}.

As for the classical CLT, the first two moments are conserved along the free
CLT: $\tau(s_n)=0$ and $\tau(s_n^2)=1$ for all $n\geq1$. The semicircle
$*$-law is the free analogue of the Gaussian distribution of classical
probability. The semicircle $*$-law is stable by free convolution: if
$a_1,\ldots,a_n$ are free with semicircle $*$-law then the $*$-law of
$a_1+\cdots+a_n$ is also semicircle and its second moment is the sum of the
second moments. In particular $s_n=(a_1+\cdots+a_n)/\sqrt{n}$ is semicircle,
just like the Gaussian case in the CLT of classical probability! If
$\mu_1,\ldots,\mu_n$ are semicircle laws then their free convolution
$\mu_1\boxplus\cdots\boxplus\mu_n$ is also a semicircle law and its variance
is the sum of the variances.

\subsection{Random walks and CLT}
Let us reconsider the free group $F_n$. The algebraic variables
\[
a_1=\frac{u_{g_1}+u_{g_1^{-1}}}{\sqrt{2}},\ldots,
a_n=\frac{u_{g_n}+u_{g_n^{-1}}}{\sqrt{2}}
\]
are free, real, centered, with unit variance. Let us define
\[
a=\sum_{i=1}^n(u_{g_i}+u_{g_i}^{-1})=\sqrt{2}\sum_{i=1}^n a_i
\quad\text{and}\quad
p
=\frac{a}{2n}
=\frac{a_1+\cdots+a_n}{\sqrt{2}n}.
\]
Then $a$ is the adjacency operator of the Cayley graph $G_n$ of the free group
$F_n$, which is $2n$-regular, without cycles (a tree!), rooted at $\varn$.
Moreover, $\DOT{p\de_v,\de_w}=\frac{1}{2n}\IND_{wv^{-1}\in S}$ where
$S=\{g_1^{\pm1},\ldots,g_n^{\pm}\}$, and thus $p$ is the transition kernel of
the simple random walk on $G_n$. For every integer $m\geq0$, the quantity
$\tau(a^m)=\DOT{a^m\de_\varn,\de_\varn}$ is the number of paths of length $m$
in $G_n$ starting and ending at the root $\varn$. From the
Kesten-McKay\footnote{Appears for regular trees of even degree (Cayley graph
  of free groups) in the doctoral thesis of Harry Kesten (1931 -- ) published
  in 1959. It seems that no connections were made at that time with the
  contemporary works \cite{MR0095527} of Eugene Wigner (1902 -- 1995) on
  random matrices published in the 1950's. In 1981, Brendan McKay (1951 -- )
  showed \cite{mckay1981expected} that these distributions appear as the
  limiting empirical spectral distributions of the adjacency matrices of
  sequences of (random) graphs which are asymptotically regular and without
  cycles (trees!). He does not cite Kesten and Wigner.} theorem
\cite{MR0109367,mckay1981expected,MR2316893}, for every integer $m\geq0$,
\[
\tau(a^m)
=\DOT{a^m\de_\varn,\de_\varn}
=\int\!x^m\,d\mu_{d}(x),
\]
where $\mu_{d}$ is the Kesten-McKay distribution with parameter $d=2n$, given
by
\[
d\mu_{d}(x)
:=\frac{d\sqrt{4(d-1)-x^2}}{2\pi(d^2-x^2)}
\mathbf{1}_{[-2\sqrt{d-1},2\sqrt{d-1}]}(x)dx.
\]
By parity we have $\tau(a^{2m+1})=0$ for every $m\geq0$. When $n=1$ then $d=2$
and $G_2$ is the Cayley graph of $F_1=\mathbb{Z}$, the corresponding
Kesten-McKay law $\mu_{2}$ is the \emph{arcsine law} on $[-2,2]$,
$$
\langle a^{2m}\delta_\varn,\delta_\varn\rangle
=\int\!x^{2m}\,d\mu_{2}(x)
=\int_{-2}^2\frac{x^{2m}}{\pi\sqrt{4-x^2}}dx
=\binom{2m}{m},
$$
and we recover the fact that the number of paths of length $2m$ in
$\mathbb{Z}$ starting and ending at the root $\varn$ (which is the
origin $0$) is given by the central binomial coefficient. The binomial
combinatorics is due to the commutativity of $F_1=\dZ$. At the opposite side
of degree, when $d=2n\to\infty$ then $\mu_{d}$, scaled by $(d-1)^{-1/2}$,
tends to the semicircle law on $[-2,2]$:
\[
\lim_{d\to\infty}
\frac{\langle a^{2m}\delta_\varn,\delta_\varn\rangle}{(d-1)^m}
=\lim_{d\to\infty}
\int\!\left(\frac{x}{\sqrt{d-1}}\right)^{2m}\!\!\!\!d\mu_{d}(x)
=\int_{-2}^2 y^{2m}\frac{\sqrt{4-y^2}}{2\pi}dy
=C_m:=\frac{1}{1+m}\binom{2m}{m}.
\]
As a consequence, we have, thanks to
$(d-1)^m\sim_{n\to\infty}d^m=(2n)^m=(\sqrt{2n})^{2m}$,
\[
\tau\PAR{\PAR{\frac{a_1+\cdots+a_n}{\sqrt{n}}}^{2m}}
=\frac{\tau(a^{2m})}{(2n)^m}
=
\frac{\DOT{a^{2m}\de_\varn,\de_\varn}}{d^m}
\underset{n\to\infty}{\longrightarrow}
C_m.
\]
This is nothing else but a free CLT for the triangular array
${((a_1,\ldots,a_n))}_{n\geq1}$! The free CLT is the algebraic structure that
emerges from the asymptotic analysis, as $n\to\infty$, of the combinatorics of
loop paths of the simple random walk on the Cayley graph $G_n$ of the free
group $F_n$, and more generally on the $d$-regular infinite graph without
cycles (a tree!) as the degree $d$ tends to infinity. 

We have $a=b_1+\cdots+b_n$ where $b_i=\sqrt{2}a_i=u_{g_i}+u_{g_i^{-1}}$. But
for any $m\in\dN$ we have $\tau(b_i)^{2m+1}=0$ and
$\tau(b_i^{2m})=\sum_{r=1}^{2m}\binom{2m}{r}\tau(u_{g_i^{2(r-m)}})=\binom{2m}{m}$,
and therefore the $*$-law of $b_i$ is the arcsine law $\mu_2$, and thanks to
the freeness of $b_1,\ldots,b_n$ we obtain
$\mu_d=\mu_{2}\boxplus\cdots\boxplus\mu_{2}$ ($d$ times).


\subsection{Free entropy}

Inspired by the Boltzmann H-Theorem view of Shannon on the CLT of classical
probability theory, one may ask if there exists, in free probability theory, a
free entropy functional, maximized by the semicircle law at fixed second
moment, and which is monotonic along the free CLT. We will see that the answer
is positive. Let us consider a real algebraic variable $a\in\cA$, $a^*=a$,
such that there exists a probability measure $\mu_a$ on $\dR$ such that for
every integer $m\geq0$,
\[
\tau(a^m)=\int\!x^m\,d\mu_a(x).
\]
Inspired from the micro-macro construction of the Boltzmann entropy, one may
consider an approximation at the level of the moments of the algebraic
variable $a$ (which is in general infinite dimensional) by Hermitian matrices
(which are finite dimensional). Namely, following Voiculescu, for every real
numbers $\veps>0$ and $R>0$ and integers $m\geq0$ and $d\geq1$, let
$\Ga_R(a;m;d,\veps)$ be the relatively compact set of Hermitian matrices
$h\in\cM_d(\dC)$ such that $\NRM{h}\leq R$ and $\max_{0\leq k\leq
  m}\ABS{\tau(a^k)-\frac{1}{d}\TR(h^k)}\leq\veps$. The volume
$\ABS{\Ga_R(a;m,d,\veps)}$ measures the degree of freedom of the approximation
of the algebraic variable $a$ by matrices, and is the analogue of the cardinal
(which was multinomial) in the combinatorial construction of the Boltzmann
entropy. We find
\[
\chi(a)
:=\sup_{R>0}\inf_{m\in\mathbb{N}}\inf_{\varepsilon>0}
\varlimsup_{d\to\infty}
\PAR{\frac{1}{d^2}\log\ABS{\Ga_R(a;m,d,\veps)}+\frac{\log(d)}{2}}
=\iint\!\log\ABS{x-y}\,d\mu_a(x)d\mu_a(y).
\]
This quantity depends only on $\mu_a$ and is also denoted $\chi(\mu_a)$. It is
a quadratic form in $\mu_a$. It is the Voiculescu entropy functional
\cite{MR1887698}. When $\mu$ is a probability measure on $\dC$, we will still
denote
\[
\chi(\mu):=\iint\!\log\ABS{x-y}\,d\mu(x)d\mu(y).
\]
However, this is not necessarily the entropy of an algebraic variable when
$\mu$ is not supported in $\dR$. The Voiculescu free entropy should not be
confused with the von Neumann entropy in quantum probability defined by
$S(a)=-\tau(a\log(a))$, which was studied by Lieb in \cite{MR0506364}. For
some models of random graphs, one can imitate Voiculescu and forge some sort
of graphical Boltzmann entropy, which can be related to a large deviations
rate functional, see \cite{bordenavecaputo} and references therein.

The semicircle law is for the free entropy the analogue of the Gaussian law
for the Boltzmann entropy. The semicircle law on $[-2,2]$ is the unique law
that maximizes the Voiculescu entropy $\chi$ among the laws on $\dR$ with
second moment equal to $1$, see \cite{AGZ}:
\[
\arg\max\BRA{\chi(\mu):\SUPP(\mu)\subset\dR,\int\!x^2\,d\mu(x)=1} %
=\frac{\sqrt{4-x^2}\IND_{[-2,2]}(x)}{2\pi}dx.
\]
How about laws on $\dC$ instead of $\dR$? The uniform law on the unit disc is
the unique law that maximizes the functional $\chi$ among the set of laws on
$\dC$ with second moment (mean squared modulus) equal to $1$, see
\cite{MR1485778} (here $z=x+iy$ and $dz=dxdy$):
\[
\arg\max\BRA{\chi(\mu):\SUPP(\mu)\subset\dC,\int\!|z|^2\,d\mu(z)=1} %
=\frac{\IND_{\{z\in\dC:|z|=1\}}}{\pi}dz,
\]
Under the uniform law on the unit disc, the real and the imaginary parts
follow the semicircle law on $[-1,1]$, and are not independent. If we say that
an algebraic variable $c\in\cA$ is circular when its $*$-law is the uniform
law on the unit disc of $\dC$, then, if $s_1,s_2\in\cA$ are free with $*$-law
equal to the semicircle law on $[-2,2]$, then $\frac{s_1+is_2}{\sqrt{2}}$ is
circular (here $i=(0,1)\in\dC$ is such that $i^2=-1$).

It turns out that the Voiculescu free entropy $\chi$ is monotonic along the
Voiculescu free CLT:
\[
\chi(a)=\chi(s_1)\leq\cdots\leq\chi(s_n)\leq\chi(s_{n+1})\leq\cdots
\underset{n\to\infty}{\nearrow}
\max\chi(s)
\]
where still $s_n=n^{-1/2}(a_1+\cdots+a_n)$ and where $s$ is an semicircle
algebraic variable. Shlyakhtenko gave a proof of this remarkable fact, based
on a free Fisher information functional (which is a Hilbert transform), that
captures simultaneously the classical and the free CLT
\cite{MR2346565,MR2304337}. The Boltzmann-Shannon H-Theorem interpretation of
the CLT is thus remarkably valid in classical probability theory, and in free
probability theory.

\subsection{A theorem of Eugene Wigner}

Let $H$ be a random $n\times n$ Hermitian matrix belonging to the Gaussian
Unitary Ensemble (GUE). This means that $H$ has density proportional to
\[
e^{-\frac{n}{4}\TR(H^2)}%
=e^{-\frac{n}{4}\sum_{1\leq j\leq 1}^n H_{jj}^2-\frac{n}{2}\sum_{1\leq j<k\leq n}|H_{jk}|^2}.
\]
The entries of $H$ are Gaussian, centered, independent, with variance $2/n$ on
the diagonal and $1/n$ outside. Let
$\mu_n=\frac{1}{n}\sum_{j=1}^n\de_{\la_{j}(H)}$ be the empirical spectral
distribution of $H$, which is a random discrete probability measure on $\dR$.
For every $m\geq0$, the mean $m$-th moment of $\mu_n$ is given by
\[
\dE\int\!x^m\,d\mu_n(x)
=\frac{1}{n}\sum_{j=1}^n\la_{j}^m
=\frac{\dE\TR(H^m)}{n}
=\sum_{i_1,\ldots,i_m}\frac{\dE(H_{i_1i_2}\cdots H_{i_{m-1}i_m}H_{i_mi_1})}{n}.
\]
In particular, the first two moments of $\dE\mu_n$ satisfy to
\[
\dE\int\!x\,d\mu_n(x)=\frac{\dE(\TR(H))}{n}
=\frac{\dE\sum_{j=1}^n H_{jj}}{n}=0
\]
and
\[
\dE\int\!x^2\,d\mu_n(x)=\frac{\dE(\TR(H^2))}{n}
=\dE\frac{\sum_{j,k=1}^n |H_{jk}|^2}{n}
=\frac{n(2/n)+(n^2-n)(1/n)}{n}
\underset{n\to\infty}{\longrightarrow}1.
\]
More generally, for any $m>2$, the computation of the limiting $m$-th moment
of $\dE\mu_n$ boils down to the combinatorics of the paths $i_1\to
i_2\to\cdots\to i_m$. It can be shown that the surviving terms as $n\to\infty$
correspond to paths forming a tree and passing exactly zero or two times per
each edge. This gives finally, for every integer $m\geq0$, denoting $C_m$ the
$m$-th Catalan number,
\[
\lim_{n\to\infty}\dE\int\!x^{2m+1}\,d\mu_n(x)=0
\quad\text{and}\quad
\lim_{n\to\infty}\dE\int\!x^{2m}\,d\mu_n(x)=C_m.
\] 
This means that $\dE\mu_n$ tends as $n\to\infty$ in the sense of moments to
the semicircle law on $[-2,2]$ (which has unit variance). Just like the CLT,
the result is actually universal, in the sense that it remains true if one
drops the Gaussian distribution assumption of the entries. This is the famous
Wigner theorem \cite{MR0095527}, in its modern general form, named after
Eugene Paul Wigner (1902 -- 1995). The GUE case is in fact exactly solvable:
one can compute the density of the eigenvalues, which turns out to be
proportional to
\[
\prod_{j=1}^ne^{-\frac{n}{4}\la_{j}^2}\prod_{j<k}(\la_j-\la_k)^2
=\exp
\PAR{-\frac{n}{4}\sum_{j=1}^n\la_j^2
-\sum_{j\neq k}\log\frac{1}{\ABS{\la_j-\la_k}}}.
\]
The logarithmic repulsion is the Coulomb repulsion in dimension $2$. Also,
this suggest to interpret the eigenvalues $\la_1,\ldots,\la_n$ as a Coulomb
gas of two-dimensional charged particles forced to stay in a one dimensional
ramp (the real line) and experiencing a confinement by a quadratic potential.
These formulas allow to deduce the semicircle limit of the one-point
correlation (density of $\dE\mu_n$), by using various methods, such as
orthogonal polynomials, or large deviations theory, see \cite{AGZ}.

\subsection{Asymptotic freeness of unitary invariant random matrices}
If $A$ and $B$ are two Hermitian $n\times n$ matrices, then the spectrum of
$A+B$ depend not only on the spectrum of $A$ and the spectrum of $B$, but also
on the eigenvectors\footnote{If $A$ and $B$ commute, then they admit the same
  eigenvectors and the spectrum of $A+B$ is the sum of the spectra of $A$ and
  $B$, but this depends on the way we label the eigenvalues, which depends in
  turn on the eigenvectors!} of $A$ and $B$. Now if $A$ and $B$ are two
independent random Hermitian matrices, there is no reason to believe that the
empirical spectral distribution $\mu_{A+B}$ of $A+B$ depend only on the
empirical spectral distributions $\mu_A$ and $\mu_B$ of $A$ and $B$. Let $A$
be a $n\times n$ random Hermitian matrix in the GUE, normalized such that
$\mu_A$ has a mean second moment equal to $1$. Then the Wigner theorem says
that $\dE\mu_A$ tend in the sense of moments, as $n\to\infty$, to the
semicircle law of unit variance. If $B$ is an independent copy of $A$, then,
thanks to the convolution of Gaussian laws, $A+B$ is identical in law to
$\sqrt{2}A$, and thus $\dE\mu_{A+B}$ tend, in the sense of moments, as
$n\to\infty$, to the semicircle law of variance $2$. Then, thanks to the free
convolution of semicircle laws, we have
\[
\dE\mu_{A+B}-\dE\mu_A\boxplus\dE\mu_B
\underset{n\to\infty}{\overset{*}{\longrightarrow}}0,
\]
where $\overset{*}{\to}$ denotes the convergence of moments. Voiculescu has
established that this \emph{asymptotic freeness} phenomenon remains actually
true beyond the GUE case, provided that the eigenspaces of the two matrices
are randomly decoupled using a random unitary conjugation. For example, let
$A$ and $B$ two $n\times n$ Hermitian matrices such that $\mu_A\to\mu_a$ and
$\mu_B\to\mu_b$ in the sense of moments as $n\to\infty$, where $\mu_a$ and
$\mu_b$ are two compactly supported laws on $\dR$. Let $U$ and $V$ be
independent random unitary matrices uniformly distributed on the unitary group
(we say Haar unitary). Then
\[
\dE\mu_{UAU^*+VBV^*}
\underset{n\to\infty}{\overset{*}{\longrightarrow}}
\mu_a\boxplus\mu_b.
\]
See \cite{AGZ}. This asymptotic freeness reveals that free probability is the
algebraic structure that emerges from asymptotic analysis of large dimensional
unitary invariant models of random matrices.

\medskip

Since the functional $\chi$ is maximized by the uniform law on the unit disc,
one may ask about an analogue of the Wigner theorem for non-Hermitian random
matrices. The answer is positive. For our purposes, we will focus on a special
ensemble of random matrices, introduced by Jean Ginibre\footnote{Jean Ginibre
  is also famous for FKG inequalities, and for scattering for Schrödinger
  operators.} in the 1960's in \cite{Ginibre}, for which one can compute
explicitly the density of the eigenvalues.

\section{Jean Ginibre and his ensemble of random matrices}

One of the most fascinating result in the asymptotic analysis of large
dimensional random matrices is the circular law for the complex Ginibre
ensemble, which can be proved using the Voiculescu functional $\chi$
(maximized at fixed second moment by uniform law on unit disc).

\subsection{Complex Ginibre ensemble}

A simple model of random matrix is the Ginibre model:
\[
G=
\begin{pmatrix}
  G_{11}&\cdots&G_{1n}\\
  \vdots & \vdots & \vdots\\
  G_{n1} & \cdots & G_{nn}
\end{pmatrix}
\]
where ${(G_{jk})}_{1\leq j,k\leq n}$ are i.i.d.\ random variables on $\dC$,
with $\Re G_{jk},\Im G_{jk}$ of Gaussian law of mean $0$ and variance
$1/(2n)$. In particular, $\dE(|G_{jk}|^2)=1/n$. The density of $G$ is
proportional to
\[
\prod_{j,k=1}^ne^{-n\ABS{G_{jk}}^2} %
= e^{-\sum_{j,k=1}^n n\ABS{G_{jk}}^2} %
= e^{-n\mathrm{Tr}(GG^*)}.
\]
This shows that the law of $G$ is unitary invariant, meaning that $UGU^*$ and
$G$ have the same law, for every unitary matrix $U$. Can we compute the law of
the eigenvalues of $G$? Let $G=UTU^*$ be the Schur unitary triangularization
of $G$. Here $T=D+N$ where $D$ is diagonal and $N$ is upper triangular with
null diagonal. In particular $N$ is nilpotent, $T$ is upper triangular, and
the diagonal of $D$ is formed with the eigenvalues $\la_1,\ldots,\la_n$ in
$\dC$ of $G$. The Jacobian of the change of variable $G\mapsto (U,D,N)$ is
proportional to $\prod_{1\leq j<k\leq n}\ABS{\la_j-\la_k}^2$ (for simplicity,
we neglect here delicate problems related to the non-uniqueness of the Schur
unitary decomposition). On the other hand,
\[
\TR(GG^*) =\TR(DD^*)+\TR(DN^*)+\TR(ND^*)+\TR(NN^*) =\TR(DD^*).
\]
This allows to integrate out the $U,N$ variables. The law of the eigenvalues
is then proportional to
\[
e^{-n\sum_{j=1}^n\ABS{\la_j}^2}\prod_{1\leq j<k\leq n}\ABS{\la_j-\la_k}^2.
\]
This defines a determinantal process on $\dC$: the complex Ginibre ensemble
\cite{MR2641363,MR2932638,MR2552864}.

\begin{figure}[htbp]
  \begin{center}
    \includegraphics[scale=.5]{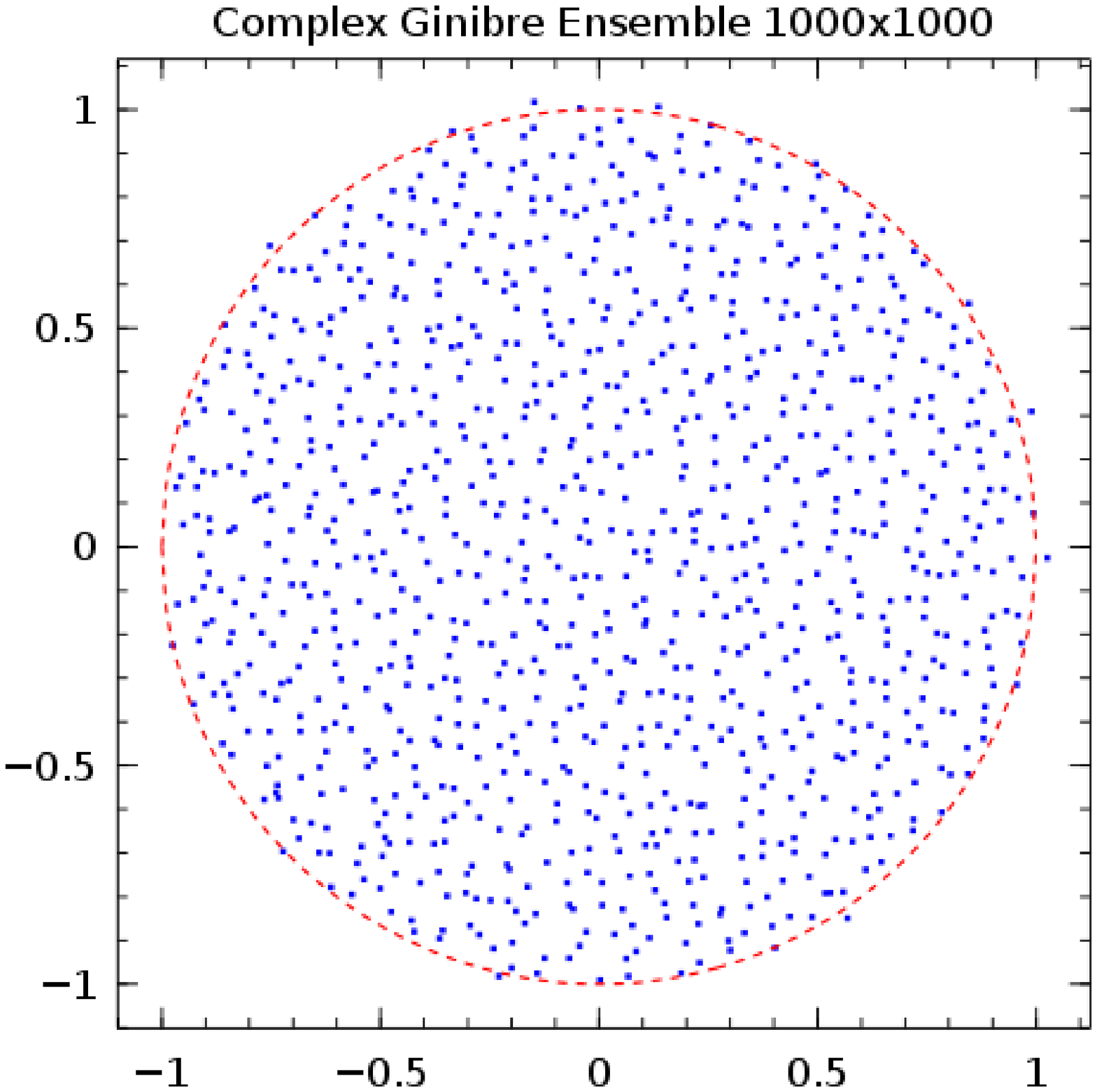}
    \caption{The eigenvalues of a single matrix drawn from the complex Ginibre
      ensemble of random matrices. The dashed line is the unit circle. This
      numerical experiment was performed using the promising Julia
      \url{http://julialang.org/}}
    \begin{lstlisting}[basicstyle=\footnotesize\ttfamily,frame=single] 
  Pkg.add("Winston"); using Winston # Pkg.add("Winston") is needed once for all.
  n=1000; (D,U)=eig((randn(n,n)+im*randn(n,n))/sqrt(2*n)); I=[-1:.01:1]; 
  J=sqrt(1-I.^2); hold(true); plot(real(D),imag(D),"b.",I,J,"r--",I,-J,"r--")
  title(@sprintf("Complex Ginibre Ensemble %dx%d",n,n)); file("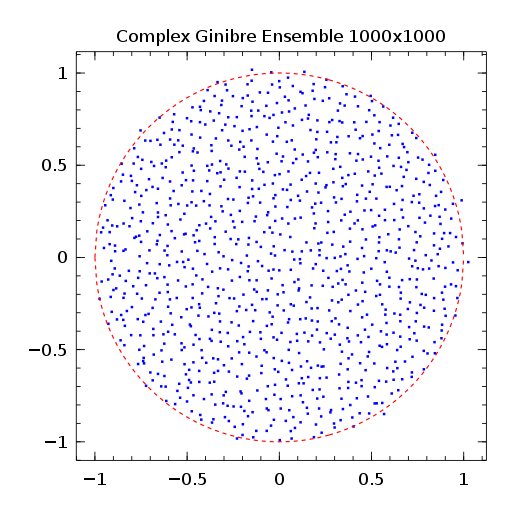")
    \end{lstlisting}
  \end{center}
\end{figure}

\subsection{Circular law for the complex Ginibre ensemble}

In order to interpret the law of the eigenvalues as a Boltzmann measure, we
put the Vandermonde determinant inside the exponential:
\[
e^{-n\sum_j\ABS{\la_j}^2+2\sum_{j<k}\log\ABS{\la_j-\la_k}}.
\]
If we encode the eigenvalues by the empirical measure
$\mu_n:=\frac{1}{n}\sum_{j=1}^n\de_{\la_j}$, this takes the form
\[
e^{-n^2 \cI(\mu_n)}
\]
where the ``energy'' $\cI(\mu_n)$ of the configuration $\mu_n$ is defined via
\[
\cI(\mu):=
\int\!\ABS{z}^2\,d\mu(z)
+\iint_{\neq}\!\log\frac{1}{\ABS{z-z'}}\,d\mu(z)d\mu(z').
\]
This suggests to interpret the eigenvalues $\la_1,\ldots,\la_n$ of $G$ as
Coulomb gas of two-dimensional charged particles, confined by a an external
field (quadratic potential) and subject to pair Coulomb repulsion. Note that
$-\cI$ can also be seen as a penalized Voiculescu functional. Minimizing a
penalized functional is equivalent to minimizing without penalty but under
constraint (Lagrange). Presently, if $\cM$ is the set of probability measures
on $\dC$ then $\inf_{\cM}\cI>-\infty$ and the infimum is achieved at a unique
probability measure $\mu_*$, which is the uniform law on the unit disc of
$\dC$.

How does the random discrete probability measure $\mu_n$ behave as
$n\to\infty$? Following Hiai and Petz \cite{MR1606719}\footnote{See also
  Anderson, Guionnet, and Zeitouni \cite{AGZ}, Ben Arous and Zeitouni
  \cite{MR1660943}, and Hardy \cite{MR2926763}.}, one may adopt a large deviations approach. Let $\cM$ be the set of
probability measures on $\dC$. One may show that the functional
$\cI:\cM\to\dR\cup\{+\infty\}$ is lower semi continuous for the topology of
narrow convergence, is strictly convex, and has compact level sets. Let us
consider a distance compatible with the topology. It can be shown that for
every ball $B$ for this distance,
\[
\dP(\mu_n\in B) \approx \exp\PAR{-n^2\inf_{B}(\cI-\inf_{\cM}\cI)}
\]
Now either $\mu_*\in B$ and in this case $\dP(\mu_n\in B)\approx 1$, or
$\mu_*\not\in B$ and in this case $\dP(\mu_n\in B)\to0$ exponentially fast.
Actually, the first Borel-Cantelli lemma allows to deduce that almost surely
\[
\lim_{n\to\infty}\mu_n 
= \mu_*
= \arg\inf \cI 
= \frac{\mathbf{1}_{\{z\in\dC:\ABS{z}\leq1\}}}{\pi}dz,
\]
where $z=x+iy$ and $dz=dxdy$. This phenomenon is known as the circular law. If
one starts with a Hermitian random Gaussian matrix -- the Gaussian Unitary
Ensemble (GUE) -- then the same analysis is available, and produces a
convergence 
to the semicircle law on $[-2,2]$.

The circular law is universal, in the sense that it remains valid if one drops
the Gaussian assumption of the entries of the matrix, while keeping the
i.i.d.\ structure and the $1/n$ variance. This was the subject of a long
series of works by Girko, Bai, Tao and Vu, among others, see
\cite{MR2906465,MR2567175,MR2908617}. Another way to go beyond the Gaussian
case is to start from the Coulomb gas and to replace the quadratic confining
potential $\ABS{\cdot}^2$ by a more general potential $V:\dC\to\dR$, not
necessarily radial. This type of generalization was studied for instance by
Saff and Totik, and by Hiai and Petz, among others, see for instance
\cite{MR1485778,MR1746976,AGZ,MR2926763}.

Beyond random matrices, how about the empirical measure of random particles in
$\dR^d$ with Coulomb type singular repulsion and external field confinement?
Is there an analogue of the circular law phenomenon? Does the ball replace
the disc? The answer is positive.

\section{Beyond random matrices}

Most of the material of this section comes from our work
\cite{2013arXiv1304.7569C} with N. Gozlan and P.-A. Zitt.

\subsection{The model}
We consider a system of $N$ particles in $\dR^d$ at positions
$x_1,\ldots,x_N$, say with charge $1/N$. These particles are subject to
confinement by an external field via a potential $x\in\dR^d\mapsto V(x)$, and
to internal pair interaction (typically repulsion) via a potential
$(x,y)\in\dR^d\times\dR^d\mapsto W(x,y)$ which is symmetric: $W(x,y)=W(y,x)$.
The idea is that an equilibrium may emerge as $N$ tends to infinity. The
configuration energy is
\begin{align*}
\cI_N(x_1,\ldots,x_N)
&=\sum_{i=1}^N\frac{1}{N}V(x_i)+\sum_{1\leq i<j\leq N}\frac{1}{N^2}W(x_i,x_j)\\
&=\int\!V(x)\,d\mu_N(x)+\frac{1}{2}\iint_{\neq}\!W(x,y)d\mu_N(x)d\mu_N(y)
\end{align*}
where $\mu_N$ is the empirical measure of the particles
(global encoding of the particle system)
\[
\mu_N:=\frac{1}{N}\sum_{k=1}^N\de_{x_k}.
\]
The model is mean-field in the sense that each particle interacts with the
others only via the empirical measure of the system. If $1\leq d\leq 2$ then
one can construct a random normal matrix which admits ours particles
$x_1,\ldots,x_N$ as eigenvalues: for any $n\times n$ unitary matrix $U$,
\[
M=U\mathrm{Diag}(x_1,\ldots,x_N)U^*,
\]
which is unitary invariant if $U$ is Haar distributed. However, we are more
interested in an arbitrarily high dimension $d$, for which no matrix model is
available. We make our particles $x_1,\ldots,x_N$ random by considering the
exchangeable probability measure $P_N$ on $(\dR^d)^N$ with density
proportional to
\[
e^{-\be_N \cI_N(x_1,\ldots,x_N)}
\]
where $\be_N>0$ is a positive parameter which may depend on $N$. The law $P_N$
is a Boltzmann measure at inverse temperature $\be_N$, and takes the form
$\prod_{i=1}^Nf_1(x_i)\prod_{1\leq i<j\leq N}f_2(x_i,x_j)$ due to the
structure and symmetries of $\cI_N$. The law $P_N$ on $(\dR^d)^N$ is informally
the invariant law\footnote{One may also view $P_N$ as the steady state of a
  Fokker-Planck evolution equation with conservation laws.} of the reversible
diffusion process ${(X_t)}_{t\in\dR_+}$ solution of the system of stochastic
differential equations
\[
dX_{t,i}=
\sqrt{\frac{2}{\be_N}}dB_{t,i}
-\frac{1}{N}\nabla V(X_{t,i})dt
-\frac{1}{N^2}\sum_{j\neq i}\nabla_{\!1}W(X_{t,i},X_{t,j})dt.
\]      
This can be seen as a special McKean-Vlasov mean-field particle system with
potentially singular interaction. The infinitesimal generator of this Markov
process is $L=\be_N^{-1}\Delta-\nabla \cI_N\cdot\nabla$. The process may
explode in finite time depending on the singularity of the interaction $W$ on
the diagonal (e.g.\ collisions of particles). The Helmholtz free energy of this
Markov process is given, for every probability density $f$ on $(\dR^d)^N$,
\[
\int_{(\dR^d)^N}\!\cI_N\,f\,dx
-\frac{1}{\be_N}\cS(f)
=\int\!\PAR{\int\!V\,d\mu_N+\frac{1}{2}\iint_{\neq}W\,d\mu_N^{\otimes 2}}\,f\,dx
-\frac{1}{\be_N}\cS(f).
\]
The model contains the complex Ginibre ensemble of random matrices as the
special case 
\[
d=2,\ \be_N=N^2,\ V(x)=\ABS{x}^2,\ W(x,y)=2\log\frac{1}{\ABS{x-y}},
\]
which is two-dimensional, with quadratic confinement, Coulomb repulsion, and
temperature $1/N^2$. Beyond this two-dimensional example, the typical
interaction potential $W$ that we may consider is the Coulomb interaction in
arbitrary dimension (we denote by $\ABS{\cdot}$ the Euclidean norm of $\dR^d$)
\[
W(x,y)=k_{\Delta}(x-y)
\quad\text{with}\quad
k_\Delta(x)
=\begin{cases}
  -\ABS{x} & \text{if $d=1$}\\
  \log\frac{1}{\ABS{x}} & \text{if $d=2$}\\
  \frac{1}{\ABS{x}^{d-2}} & \text{if $d\geq3$}
\end{cases}
\]
and the Riesz interaction, $0<\alpha<d$ (Coulomb if $d\geq3$ and $\al=2$) $d\geq1$
\[
W(x,y)=k_{\Delta_\alpha}(x-y) 
\quad\text{with}\quad
k_{\Delta_\alpha}(x)=\frac{1}{\ABS{x}^{d-\alpha}}.
\]    
The Coulomb kernel $k_\De$ is the fundamental solution of the Laplace
equation, while the Riesz kernel $k_{\De_\al}$ is the fundamental solution of
the fractional Laplace equation, hence the notations. In other words, in the
sense of Schwartz distributions, for some constant $c_d$,
\[
\De_\al k_{\De_\al}=c_d\de_0.
\]
If $\al\neq2$ then the operator $\De_\al$ is a non-local Fourier multiplier.
 

\subsection{Physical control problem}

With the Coulomb-Gauss theory of electrostatic phenomena in mind, it is
tempting to consider the following physical control problem: given an internal
interaction potential $W$ and a target probability measure $\mu_*$ in $\dR^d$,
can we tune the external potential $V$ and the cooling scheme $\be_N$ is order
to force the empirical measure $\mu_N=\frac{1}{N}\sum_{i=1}^N\de_{x_i}$ to
converge to $\mu_*$ as $N\to\infty$? One can also instead fix $V$ and seek for
$W$.

\subsection{Topological result: large deviations principle}

The confinement is always needed in order to produce a non degenerate
equilibrium. This is done here by an external field. This can also be done by
forcing a compact support, as Frostman did in his doctoral thesis
\cite{frostman}, or by using a manifold as in Dragnev and Saff
\cite{MR2276529} and Berman \cite{2008arXiv0812.4224B}. Let $\cM_1$ be the set
of probability measures on $\dR^d$ equipped with the topology of narrow
convergence, which is the dual convergence with respect to continuous and
bounded test functions.
From the expression of $\cI_N$ in terms of $\mu_N$, the natural limiting energy
functional is the quadratic form with values in $\dR\cup\{+\infty\}$ defined by
\[
\mu\in\cM_1\mapsto
\cI(\mu)=\int\!V(x)\,d\mu(x)+\frac{1}{2}\iint\!W(x,y)\,d\mu(x)d\mu(y).
\]
We make the following assumptions on $V$ and $W$ (fulfilled for instance when
the localization potential is quadratic $V=c\ABS{\cdot}^2$ and the interaction
potential $W$ is Coulomb or Riesz).
\begin{itemize}
\item \emph{Localization and repulsion.}
  \begin{itemize}
  \item the function $V:\dR^d\to\dR$ continuous, $V(x)\to+\infty$ as
    $\ABS{x}\to+\infty$, $e^{-V}\in L^1(dx)$;
  \item the function $W:\dR^d\times\dR^d\to\dR\cup\{+\infty\}$ continuous,
    finite outside diagonal, and symmetric $W(x,y)=W(y,x)$, (however $W$ can
    be infinite on the diagonal!);
  \end{itemize}
\item \emph{Near infinity confinement beats repulsion.} For some constants
  $c\in\dR$ and $\veps_o\in(0,1)$,
  \[
  \forall x,y\in\dR^d,\quad W(x,y)\geq c-\veps_o(V(x)+V(y));
  \]
\item \emph{Near diagonal repulsion is integrable.} for every compact
  $K\subset\dR^d$,
  \[
  z\mapsto\sup_{x,y\in K,\ABS{x-y}\geq\ABS{z}}W(x,y)\in L^1(dz);
  \]
\item \emph{Regularity.} $\forall\nu\in\cM_1(\dR^d)$, if $\cI(\nu)<\infty$
  then
  \[
  \exists(\nu_n)\in\cM_1(\dR^d),
  \quad\nu_n\ll dx,\quad \nu_n\to\nu, \quad \cI(\nu_n)\to \cI(\nu);
  \]
\item \emph{Cooling scheme.} $\be_N\gg N\log(N)$ (for the Ginibre ensemble
  $\be_N=N^2$).
\end{itemize}
Under these assumptions, it is proven in \cite{2013arXiv1304.7569C} that the
sequence ${(\mu_N)}_{N\geq1}$ of random variables taking values in $\cM_1$
satisfies to a Large Deviations Principle (LDP) at speed $\be_N$ with good
rate function $\cI-\inf_{\cM_1}\cI$. In other words:


\begin{itemize}
\item Rate function: $\cI$ is lower semi-continuous with compact level sets, and
  $\inf_{\cM_1}\cI>-\infty$;
\item LDP lower and upper bound: for every Borel set $A$ in $\cM_1$,
  \[
  \liminf_{N\to\infty}
  \frac{\log P_N(\mu_N\in A)}{\be_N}
  \geq-\inf_{\mu\in\mathrm{int}(A)}(\cI-\inf \cI)(\mu)
  \]
  and
  \[
  \limsup_{N\to\infty}
  \frac{\log P_N(\mu_N\in A)}{\be_N}
  \leq-\inf_{\mu\in\mathrm{clo}(A)}(\cI-\inf \cI)(\mu);
  \]    
\item Convergence: $\arg\inf \cI$ is not empty and almost surely
  $\lim_{N\to\infty}\mathrm{dist}(\mu_N,\arg\inf \cI)=0$ where $\mathrm{dist}$
  is the bounded-Lipschitz dual distance (it induces the narrow topology on
  $\cM_1$).
\end{itemize}

This LDP must be seen as an attractive tool in order to show the convergence
of $\mu_N$. This topological result is built on the idea that the density of
$P_N$ is proportional as $N\gg1$ to
\[
e^{-\be_N \cI_N(\mu_N)} \approx e^{-\be_N\cI(\mu_N)},
\]     
and thus, informally, the first order global asymptotics as $N\gg1$ is
\[
\mu_N\approx\arg\inf \cI.
\]
This generalizes the case of the complex Ginibre ensemble considered in the
preceding section. At this level of generality, the set $\arg\inf \cI$ is
non-empty but is not necessarily a singleton. In the sequel, we provide a more
rigid differential result in the case where $W$ is the Riesz potential, which
ensures that $\cI$ admits a unique minimizer, which is characterized by simple
properties.

The $N\log(N)$ in condition $\be_N\gg N\log(N)$ comes from volumetric
(combinatorial) estimates.

It is important to realize that $W$ is potentially singular on the diagonal,
and that this forbids the usage of certain LDP tools such as the
Laplace-Varadhan lemma or the Gärtner-Ellis theorem, see \cite{MR2571413}. If
$W$ is continuous and bounded on $\dR^d\times\dR^d$ then one may deduce our
LDP from the LDP with $W\equiv0$ by using the Laplace-Varadhan lemma.
Moreover, if $W\equiv0$ then $P_N=\eta_N^{\otimes N}$ is a product measure
with $\eta_N\propto e^{-(\be_N/N)V}$, the particles are independent, the rate
functional is
\[
\cI(\mu)-\inf_{\cM_1}\cI=\int\!V\,d\mu-\inf V,
\quad\text{and}\quad
\arg\inf_{\cM_1} \cI %
=\cM_V =\{\mu\in\cM_1:\mathrm{supp}(\mu)\subset\arg\inf V\},
\]
which gives, thanks to $\be_N\gg N\log(N)\gg N$,
\[
\lim_{N\to\infty}\mathrm{dist}(\mu_N,\cM_V)=0.
\]

\subsection{Linear inverse temperature and link with Sanov theorem}

If $\be_N=N$ and $W\equiv0$ then $P_N=(\mu_*)^{\otimes N}$ is a product
measure, the law $\mu_*$ has density proportional to $e^{-V}$, the particles
are i.i.d.\ of law $\mu_*$, and the Sanov theorem states that
${(\mu_N)}_{N\geq1}$ satisfies to an LDP in $\cM_1$ with rate function given
by the Kullback-Leibler relative entropy $\cK$ with respect to $\mu_*$.
If $W\in\cC_b$ then the particles are no longer independent but the
Laplace-Varadhan lemma allows to deduce that ${(\mu_N)}_{N\geq1}$ satisfies to
an LDP in $\cM_1$ with rate function $\cR$ given by
\[
\cR(\mu) %
= \cK(\mu) + \frac{1}{2}\iint\!W(x,y)\,d\mu(x)d\mu(y) 
= -\cS(\mu) + \cI(\mu)
\]
where $\cS(\mu)=\cS(f)$ if $d\mu=fd\mu_*$ and $\cS(\mu)=+\infty$ otherwise. We
have here a contribution of the Boltzmann entropy $\cS$ and a Voiculescu type
functional $\chi$ via its penalized version $\cI$. Various versions of this
LDP was considered in the literature in various fields and in special
situations, for instance in the works of Messer and Spohn \cite{MR704588},
Kiessling \cite{MR1193342}, Caglioti, Lions, Marchioro, and Pulvirenti
\cite{MR1145596}, Bodineau and Guionnet \cite{MR1678526}, and Kiessling and
Spohn \cite{MR1669669}, among others.

\subsection{Differential result: rate function analysis}

Recall that $\cI$ is a quadratic form on measures. 
\[
\frac{t\cI(\mu)+(1-t)\cI(\nu)-\cI(t\mu+(1-t)\nu)}{t(1-t)}
=\iint\!W\,d(\mu-\nu)^2.
\]
This shows that $\cI$ is convex if and only if $W$ is weakly positive in the
sense of Bochner. Term ``weak'' comes from the fact that we need only to check
positivity on measures which are differences of probability measures. The
Coulomb kernel in dimension $d=2$ is not positive, but is weakly positive. It
turns out that every Coulomb or Riesz kernel is weakly positive, in any
dimension. Consequently, the functional $\cI$ is convex in the case of Coulomb
or Riesz interaction, for every dimension $d$. One may also rewrite the
quadratic form $\cI$ as
\[
\cI(\mu)=\int\!V\,d\mu+\frac{1}{2}\int\!U_\mu\,d\mu
\quad\text{where}\quad
U_\mu(x):=\int\!W(x,y)\,d\mu(y).
\]
In the Coulomb case then $U_\mu$ is the logarithmic potential in dimension
$d=2$, and the Coulomb or minus the Newton potential in higher dimensions. An
infinite dimensional Lagrange variational analysis gives that the gradient of
the functional $\cI$ at point $\mu$ is $V+U_\mu$ should be constant on the
support of the optimum, and in the case where $W(x,y)=k_D(x-y)$ where $k_D$ is
the fundamental solution of a local (say differential) operator $D$, meaning
$Dk_D=-\de_0$ in the sense of distributions, we have $DU_\mu=-\mu$ which gives
finally that on the support of $\mu_*$,
\[
\mu_*=DV,
\]
In particular, if $V$ is the squared Euclidean norm and if $W$ is the Coulomb
repulsion then $D$ is the Laplacian and $DV$ is constant, which suggests that
$\mu_*$ is constant on its support (this is compatible with what we already
know for dimension $d=2$, namely $\mu_*$ is uniform on the unit disc of
$\dC$).

We show in \cite{2013arXiv1304.7569C} that if $W$ is the Riesz interaction,
then the functional $\cI$ is strictly convex, and $\arg\inf \cI=\{\mu_*\}$ is a
singleton, and $\mu_*$ is compactly supported, and almost surely,
\[
\mu_N\underset{N\to\infty}{\longrightarrow}\mu_*,
\]
and moreover $\mu_*$ is characterized by the Lagrange conditions
\[
U_{\mu_*}+V=C_*\text{ on supp$(\mu_*)$ and $\geq C_*$ outside}
\]
In the Coulomb case the constant $C_*$ is known as the modified Robin
constant. The support constraint in the Lagrange conditions make difficult the
analysis of the Riesz case beyond the Coulomb case, due to the non-local
nature of the fractional Laplacian. Finally, let us mention that it is shown
in \cite{2013arXiv1304.7569C} using the Lagrange conditions that one can
construct $V$ from $\mu_*$ if $\mu_*$ is smooth enough, and this gives a
positive answer to the physical control problem mentioned before.

In the Coulomb case, we have $\De U_\mu=c_d\mu$, and an integration parts
gives
\[
\cI(\mu)-\int\!V\,d\mu
=\frac{1}{2}\int\!U_\mu\,d\mu
=\frac{1}{2}\int\!U_\mu\,\De U_\mu\,dx
=\frac{1}{2}\int\!\ABS{\nabla U_\mu}^2\,dx.
\]
The right hand side is the integral of the squared norm of the gradient
$\nabla U_\mu$ of the electrostatic potential $U_\mu$ generated by $\mu$, in
other words the ``squared-field'' (``carré du champ'' in French).


\begin{table}
  \begin{tabular}{r|l}
    $\cdots$ & $\cdots$\\
     1660 & Newton \\
     1760 & Coulomb, Euler \\
     1800 & Gauss \\
     1820 & Carnot \\
     1850 & Helmholtz \\
     1860 & Clausius\\
     1870 & Boltzmann, Gibbs, Maxwell\\
     1900 & Markov, Perron and Frobenius \\
     1930 & Fisher, Kolmogorov, Vlasov \\
     1940 & de Bruijn, Kac, von Neumann, Shannon \\
     1950 & Linnik, Kesten, Kullback, Sanov, Stam, Wigner\\
     1960 & Ginibre, McKean, Nelson \\
     1970 & Cercignani, Girko, Gross \\
     1980 & Bakry and Émery, McKay, Lions, Varadhan, Voiculescu\\
     2000 & Perelman, Tao and Vu, Villani \\
     $\cdots$ & $\cdots$
  \end{tabular}
  \caption{The arrow of time and some of the main actors mentioned in the
    text. One may have in mind the Stigler law: ``\emph{No scientific
      discovery is named after its original discoverer.}'', which is attributed
    by Stephen Stigler to Robert K. Merton.}
\end{table}

\subsection{Related problems}

In the Coulomb case, and when $V$ is radial, then $\mu_*$ is supported in a
ring and one can compute its density explicitly thanks to the Gauss averaging
principle, which states that the electrostatic potential generated by a
distribution of charges in a compact set is, outside the compact set, equal to
the electrostatic potential generated by a unique charge at the origin, see
\cite{MR2647570}. In particular, in the Coulomb case, and when $V$ is the
squared norm, then $\mu_*$ is the uniform law on a ball of $\dR^d$,
generalizing the circular law of random matrices.

Beyond the Coulomb case, even in the Riesz case, no Gauss averaging principle
is available, and no one knows how to compute $\mu_*$. Even in the Coulomb
case, the computation of $\mu_*$ is a problem if $V$ is not rotationally
invariant. See Saff and Dragnev \cite{MR2276529}, and Bleher and Kuijlaars
\cite{MR2921180}.

When $V$ is weakly confining, then the LDP may be still available, but $\mu_*$
is no longer compactly supported. This was checked in dimension $d=2$ with the
Coulomb potential by Hardy in \cite{MR2926763}, by using a compactification
(stereographic) which can be probably used in arbitrary dimensions.

It is quite natural to ask about algorithms (exact or approximate) to
simulated the law $P_N$. The Coulomb case in dimension $d=2$ is determinantal
and this rigid structure allows exact algorithms \cite{MR2216966}. Beyond this
special situation, one may run an Euler-Langevin MCMC approximate algorithm
using the McKean-Vlasov system of particles. How to do it smartly? Can we do
better? Speaking about the McKean-Vlasov system of particles, one may ask
about its behavior when $t$ and/or $N$ are large. How it depends on the
initial condition, at which speed it converges, do we have a propagation of
chaos, does the empirical measure converge to the expected PDE, is it a
problem to have singular repulsion? Can we imagine kinetic versions in
connection with recent works on Vlasov-Poisson equations for instance
\cite{MR2065020}? Similar questions are still open for models with attractive
interaction such as the Keller-Segel model.

Beyond the first order global analysis, one may ask about the behavior of
$\mu_N-\mu_*$. This may lead to central limit theorems in which the speed may
depend on the regularity of the test function. Some answers are already
available in dimension $d=2$ in the Coulomb case by Ameur, Hedenmalm, and
Makarov \cite{MR2817648}. Another type of second order analysis is studied by
Serfaty and her collaborators \cite{2014arXiv1403.6860S,MR3163544}.
Similar infinite dimensional mean-field models are studied by Lewin and
Rougerie in relation with the Hartree model for Bose-Einstein condensates
\cite{lewin-gazette}.

Such interacting particle systems with Coulomb repulsion have inspired similar
models in complex geometry, in which the Laplacian is replaced by the
Monge-Ampere equation. The large deviations approach remains efficient in this
context, and was developed by Berman \cite{2008arXiv0812.4224B}.

The limit in the first order global asymptotics depend on $V$ and $W$, and is
thus non universal. In the case of $\beta$-ensembles of random matrix theory
(particles confined on the real line $\dR$ with Coulomb repulsion of dimension
$2$), the asymptotics of the local statistics are universal, and this was the
subject of many contributions, with for instance the works of Ramirez, Rider,
and Virág \cite{MR2813333}, Erdős, Schlein, and Yau \cite{MR2810797},
Bourgade, Erdős, and Yau \cite{MR2905803}, and Bekerman, Figalli, and Guionnet
\cite{BFG}. Little is known in higher dimensions or with other interaction
potential $W$. The particle of largest norm has Tracy-Widom fluctuation in
$\be$-ensembles, and here again, little is known in higher dimensions or with
other interactions, see \cite{chapec} for instance and references therein.


\subsection{Caveat} 

The bibliography provided in this text is informative but incomplete. We
emphasize that this bibliography should not be taken as an exhaustive
reference on the history of the subject.

\subsection{Acknowledgments} 

The author would like to thank Laurent Miclo and Persi Diaconis for their
invitation, and François Bolley for his useful remarks on an early draft of
the text. This version benefited from the constructive comments of two
anonymous reviewers.

\begin{center}
  \begin{figure}[htbp]
    \includegraphics[scale=0.4]{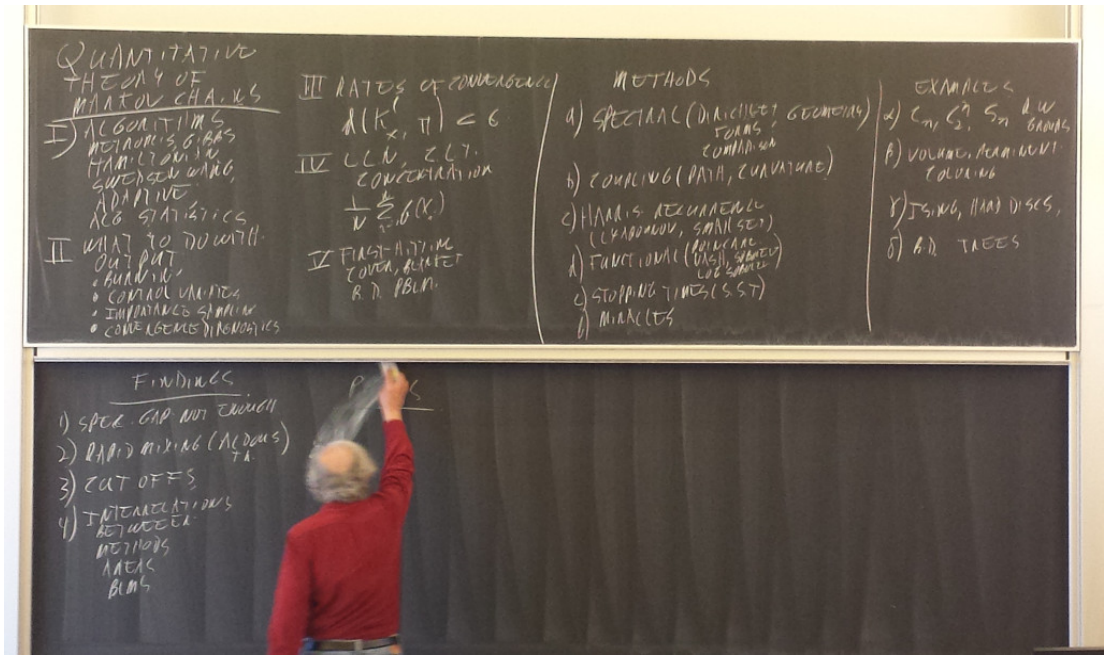}
    \caption*{Persi Diaconis during his talk. Institut de Mathématiques de Toulouse, March 24, 2014.}
  \end{figure}
\end{center}

\makeatletter
\def\@MRExtract#1 #2!{#1}     
\renewcommand{\MR}[1]{
  \xdef\@MRSTRIP{\@MRExtract#1 !}%
  \href{http://www.ams.org/mathscinet-getitem?mr=\@MRSTRIP}{MR-\@MRSTRIP}}
\makeatother

\end{document}